\documentclass[11 pt]{amsart}
\usepackage[pagebackref]{hyperref}
\hypersetup{
colorlinks = true,
linkcolor = blue,
citecolor = blue}
\usepackage[utf8]{inputenc}
\usepackage[T1]{fontenc}

\usepackage[]{amsfonts}
\usepackage{amssymb}
\usepackage{amsmath}
\def\couleur(#1 #2 #3)
	{
\special{ps::[end]}}

\def\bx#1{\setbox1=\hbox{\kern3pt{#1}\kern3pt}			% Make a box. Close it by "}"
 \dimen1=\ht1 \advance\dimen1 by 3pt \dimen2=\dp1 \advance\dimen2 by 3pt
 \setbox1=\hbox{\vrule height\dimen1 depth\dimen2\box1\vrule}%
 \setbox1=\vbox{\hrule\box1\hrule}%
 \advance\dimen1 by .4pt \ht1=\dimen1
 \advance\dimen2 by .4pt \dp1=\dimen2 \box1\relax}

\def\wbb#1{\kern#1em}
\def\vci{\vrule  width.02em height1.47ex depth-.0ex}		% le 1 en blackboard
\def\11{{\rm\wbb{.2}\vci\wbb{-.37}1}}

\def\underset#1#2{\mathrel{\mathop{\kern0pt #2}\limits_{#1}}}

\def\overset#1#2{\mathrel{\mathop{\kern0pt #2}\limits^{#1}}}

\newtheorem{thm}{Theorem}[section]
\newtheorem{lem}[thm]{Lemma}
\newtheorem{prop}[thm]{Proposition}
\newtheorem{cor}[thm]{Corollary}
\newtheorem{defin}[thm]{Definition}
\newtheorem{rem}[thm]{Remark}
\newtheorem{exa}[thm]{Example}

\begin{document}

\title{The LIR Method. $L^{r}$ solutions of elliptic equation in a complete riemannian manifold.}

\author{Eric Amar}

\date{}

\subjclass[2010]{35J58 58J05 58A14}

\keywords{Elliptic linear equation, riemannian manifold, Sobolev estimates}
\maketitle
 \ \par 
\ \par 

\tableofcontents
\ \par 
\renewcommand{\abstractname}{Abstract}

\begin{abstract}
\quad We introduce the Local Increasing Regularity Method (LIRM) which
 allows us to get from \emph{local} a priori estimates, on solutions
 $u$ of a linear equation $\displaystyle Du=\omega ,$ \emph{global} ones.\ \par 
\quad As an application we shall prove that if $D$ is an elliptic linear
 differential operator of order $m$ with ${\mathcal{C}}^{\infty
 }$ coefficients operating on the sections of a complex vector
 bundle $\displaystyle G:=(H,\pi ,M)$ over a compact Riemannian
 manifold $M$ without boundary and $\omega \in L^{r}_{G}(M)\cap
 (\mathrm{k}\mathrm{e}\mathrm{r}D^{*})^{\perp },$ then there
 is a $u\in W^{m,r}_{G}(M)$ such that $Du=\omega $ on $M.$\ \par 
\quad Next we investigate the case of a compact manifold with boundary
 by use of the "riemannian double manifold".\ \par 
\quad In the last sections we study the more delicate case of a complete
 but non compact Riemannian manifold by use of adapted weights.\ \par 
\end{abstract}
\ \par 
\ \par 
\ \par 
\ \par 

\section{Introduction.}
\quad Let $(M,g)$ be a complete Riemannian manifold and $\Delta :=dd^{*}+d^{*}d$
 be the Hodge laplacian on it. Let $\displaystyle \Lambda ^{p}(M)$
 be the set of $p$-forms ${\mathcal{C}}^{\infty }$ smooth on
 $M,$ then we have $\Delta :\ \Lambda ^{p}\rightarrow \Lambda
 ^{p}.$ The Poisson equation $\Delta u=\omega $ for $\omega \in
 \Lambda ^{p}(M)$ was extensively studied. Set $L_{p}^{r}$ the
 closure of $\displaystyle \Lambda ^{p}(M)$ in the space $L^{r}(M)$
 for the volume measure of $M.$ We define as usual the Sobolev
 spaces $W_{p}^{k,r}(M)$ to be the set of $p$-forms on $M$ in
 $L^{r}_{p}(M)$ together with all its covariant derivatives up
 to order $k.$\ \par 
Then $\displaystyle L_{p}^{r}$ estimates for the solutions of
 the Poisson equation are essentially equivalent to the $L_{p}^{r}$
 Hodge decomposition:\ \par 
\ \par 
\quad \quad \quad $\displaystyle L^{r}_{p}(M)={\mathcal{H}}^{r}_{p}\oplus dW^{1,r}_{p-1}(M)\oplus
 d^{*}W^{1,r}_{p+1}(M).$\ \par 
\ \par 
\quad Let us recall some results in the case $M$ compact without boundary.\ \par 
The basic work of CB Morrey~\cite{Morrey66} for $\omega \in L^{2}(M)$
 has lead to the $L^{2}$ Hodge decomposition:\ \par 
\ \par 
\quad \quad \quad $\displaystyle L^{2}_{p}(M)={\mathcal{H}}^{2}_{p}\oplus dW^{1,2}_{p-1}(M)\oplus
 d^{*}W^{1,2}_{p+1}(M)$\ \par 
\ \par 
useful in Algebraic Geometry, see C. Voisin~\cite{Voisin02}.\ \par 
\quad In 1995  Scott~\cite{Scott95} proved a strong $\displaystyle
 L^{r}$ Hodge decomposition:\ \par 
\ \par 
\quad \quad \quad $\displaystyle \forall r>1,\ L^{r}_{p}(M)={\mathcal{H}}^{r}_{p}\oplus
 dW^{1,r}_{p-1}(M)\oplus d^{*}W^{1,r}_{p+1}(M).$\ \par 
\ \par 
G. Schwarz~\cite{Schwarz95} proved the same result but in a compact
 riemannian manifold with boundary.\ \par 
\ \par 
For the case of a complete non compact riemannian manifold there
 are also classical results.\ \par 
\quad In 1949, Kodaira~\cite{Kodaira49} proved that the $\displaystyle
 L^{2}$-space of $p$-forms on $\displaystyle (M,g)$ has the (weak)
 orthogonal decomposition:\ \par 
\ \par 
\quad \quad \quad $L^{2}_{p}(M)={\mathcal{H}}^{2}_{p}\oplus {\overline{d{\mathcal{D}}_{p-1}(M)}}\oplus
 {\overline{d^{*}{\mathcal{D}}_{p+1}(M)}},$\ \par 
\ \par 
and in 1991 Gromov~\cite{Gromov91} proved a strong $\displaystyle
 L^{2}$ Hodge decomposition, under the hypothesis that $\Delta
 $ has a spectral gap in $\displaystyle L_{p}^{2}$:\ \par 
\ \par 
\quad \quad \quad $\displaystyle L^{2}_{p}(M)={\mathcal{H}}^{2}_{p}\oplus dW^{1,2}_{p-1}(M)\oplus
 d^{*}W^{1,2}_{p+1}(M).$\ \par 
\ \par 
There are also nice results by X-D. Li~\cite{XDLi09} who proved
 a strong $\displaystyle L^{r}$ Hodge decomposition on complete
 non compact riemannian manifold. See the references list on
 these questions therein.\ \par 
Finally, by use of the raising steps method, I proved in~\cite{HodgeNonCompact15},
 that we have a \emph{non classical weighted} $L^{r}_{p}(M)$
 Hodge decomposition in a complete non compact riemannian manifold.\ \par 
\ \par 
\quad The aim of this work is to extend these results to the general
 case of a linear elliptic operator $D$ of order $m$ in place
 of the Hodge laplacian. If $(M,g)$ is a compact boundary-less
 riemannian manifold, this was done in the $L^{2}$ case, for
 instance, by Warner~\cite{Warner83} and Donaldson~\cite{Donaldson08}.
 See the references therein.\ \par 
\quad Here we shall study the equation $Du=\omega $ for a general linear
 elliptic operator $D$ of order $m$ acting on sections of  $\displaystyle
 G:=(H,\pi ,M),$ a complex ${\mathcal{C}}^{m}$ vector bundle
 over $M$ of rank $N$ with fiber $H$  in the riemannian manifold $M.$\ \par 
\quad Let $M$ be a complete $n$-dimensional ${\mathcal{C}}^{m}$ riemannian
 manifold for some $m\in {\mathbb{N}},$  and let $G:=(H,\pi ,M)$
 be a complex ${\mathcal{C}}^{m}$ vector bundle over $M$ of rank
 $N$ with fiber $H.$ By a trivializing coordinate system $(U_{\varphi
 },\varphi ,\ \chi _{\varphi })$ for $G$ we mean a chart $\varphi
 $ of $M$ with domain $U_{\varphi }\subset M$ together with a
 trivializing map:\ \par 
\quad \quad \quad $\pi ^{-1}(U_{\varphi })\rightarrow U_{\varphi }{\times}H,\ g\rightarrow
 (\pi (g),\chi _{\varphi }(g))$\ \par 
over $U_{\varphi }$ for $G.$ Given a section $u$ of $G,$ its
 local representation $u_{\varphi }$ with respect to $\displaystyle
 (U_{\varphi },\varphi ,\ \chi _{\varphi })$ is defined by  
       $u_{\varphi }:=\chi _{\varphi }\circ u\circ \varphi ^{-1}.$\ \par 
Then given $s\in \lbrack 0,\ m\rbrack $ and $r\in (1,\infty ),$
 we denote by $W^{s,r}_{G}(M)$ the vector space of all sections
 $u$ of $G$ such that $\psi u_{\varphi }\in W^{s,r}(\varphi (U_{\varphi
 }),H)$ for each ${\mathcal{C}}^{m}$ function $\psi $ with compact
 support in $\varphi (U_{\varphi })\subset {\mathbb{R}}^{n}$
 and each trivializing coordinate system $(\varphi ,U_{\varphi
 },\chi _{\varphi })$ for $G,$ where sections coinciding almost
 everywhere have been identified and $\displaystyle W^{s,r}$
 is the usual Sobolev space whose main properties are recalled
 in the Subsection~\ref{CL27} of the Appendix.  In particular
 we have $\displaystyle L^{r}_{G}(M)=W^{0,r}_{G}(M).$ \ \par 
By analogy with the bundle of $p$-forms on $M,$ we shall call
 $G$-forms  the measurable sections of $G.$\ \par 
\ \par 
\quad The method we shall use is different from the previous ones.
 We shall provide a way to go from \emph{local results} to \emph{global
 ones} by use of the Local Increasing Regularity, LIR for short,
 given by the fundamental elliptic estimates. We shall introduce
 a quite general method, the LIR method, which allows us to get
 the generalization to $L^{r}$ of the result of Warner~\cite{Warner83}
 and Donaldson~\cite{Donaldson08} done for $L^{2}.$\ \par 

\begin{thm}
Let $(M,g)$ be a ${\mathcal{C}}^{\infty }$ smooth compact riemannian
 manifold without boundary. Let $D:G\rightarrow G$ be an elliptic
 linear differential operator of order $m$ with ${\mathcal{C}}^{\infty
 }$ coefficients acting on the complex ${\mathcal{C}}^{m}$ vector
 bundle $G$ over $M.$  Let $\omega \in L^{r}_{G}(M)\cap \mathrm{(}\mathrm{k}\mathrm{e}\mathrm{r}D^{*}\mathrm{)}^{\perp
 }$ with $r\geq 2.$ Then there is a bounded linear operator $S:L^{r}_{G}(M)\cap
 \mathrm{(}\mathrm{k}\mathrm{e}\mathrm{r}D^{*}\mathrm{)}^{\perp
 }\rightarrow W^{m,r}_{G}(M)$ such that $\displaystyle DS(\omega
 )=\omega $ on $M.$ So, with $u:=S\omega $ we get $\displaystyle
 Du=\omega $ and $\displaystyle u\in W^{m,r}_{G}(M).$
\end{thm}
\quad By duality we get the range $r<2$ as we did in~\cite{AmarSt13},
 using an avatar of the Serre duality~\cite{Serre55}.\ \par 
\ \par 
\quad To study the same problem when $M$ has a smooth boundary $\partial
 M,$ we shall use the technique  of the "Riemannian double".\ \par 
\quad The "Riemannian double" $\Gamma :=\Gamma (M)$ of $M,$ obtained
 by gluing two copies of (a slight extension of) $M$ along $\partial
 M,$ is a compact Riemannian manifold without boundary. Moreover,
 by its very construction, it is always possible to assume that
 $\displaystyle \Gamma $ contains an isometric copy $M$ of the
 original domain $M.$ See Guneysu and Pigola ~\cite[Appendix
 B]{GuneysuPigola}.\ \par 
\quad We shall need:\ \par 

\begin{defin}
~\label{LIR38}We shall say that $D$ has the weak maximum property,
 WMP, if, for any smooth $DG$-harmonic $h,$ i.e. a $G$-form such
 that $Dh=0$ in $M,$ smooth up to the boundary $\partial M,$
 which is flat on $\partial M,$ i.e. zero on $\partial M$ with
 all its derivatives, then $h$ is zero in $M.$
\end{defin}
\quad This definition has to be linked to the Definition~\cite[Introduction,
 p. 948]{Kenig86}:\ \par 

\begin{defin}
We shall say that an operator $D$ has the Unique Continuation
 Property, UCP, if $Du=0$ on $\Gamma $ and $u=0$ in an open set
 ${\mathcal{O}}\neq \emptyset $ of $\Gamma $ implies that $u\equiv
 0$ in $\Gamma .$
\end{defin}
\quad WMP is weaker than the UCP, because if $D$ has the UCP and if
 $h$ is flat on $\partial M,$ then we can extend $h$ by zero
 in $M^{c}$ in $\Gamma ,$ which makes $h$ still $DG$-harmonic,
 and apply the UCP to get that $h$ is zero in $M.$\ \par 
\quad The Hodge laplacian in a riemannian manifold has the UCP for
 $p$-forms by a difficult result  by N. Aronszajn, A. Krzywicki
 and J. Szarski~\cite{Aronszajn62}. Then we get:\ \par 

\begin{thm}
Let $M$ be a smooth compact riemannian manifold with smooth boundary
 $\partial M.$ Let $\displaystyle \omega \in L^{r}_{G}(M).$ There
 is a form $u\in W^{m,r}_{G}(M),$ such that $Du=\omega $ and
 ${\left\Vert{u}\right\Vert}_{W^{m,r}_{G}(M)}\leq c{\left\Vert{\omega
 }\right\Vert}_{L^{r}_{G}(M)},$ provided that the operator $D$ has the WMP.
\end{thm}
\ \par 
\quad We shall use the same ideas as we did in~\cite{HodgeNonCompact15}
 to go from the compact case to the non compact one.\ \par 
First we have to define a $m,\epsilon $-admissible ball centered
 at $x\in M.$ Its radius $R(x)$ must be small enough to make
 that ball like its euclidean image. Precisely:\ \par 

\begin{defin}
Let~\label{mLIR25} $(M,g)$ be a riemannian manifold and $\displaystyle
 x\in M.$ We shall say that the geodesic ball $\displaystyle
 B(x,R)$ is $m,\epsilon $ {\bf admissible} if there is a chart
 $\varphi \ :\ (y_{1},...,y_{n})\rightarrow {\mathbb{R}}^{n}$
 defined on it with\par 
\quad 1) $\displaystyle (1-\epsilon )\delta _{ij}\leq g_{ij}\leq (1+\epsilon
 )\delta _{ij}$ in $\displaystyle B(x,R)$ as bilinear forms,\par 
\quad 2) $\displaystyle \ \sum_{\left\vert{\beta }\right\vert \leq
 m-1}{\sup \ _{i,j=1,...,n,\ y\in B_{x}(R)}\left\vert{\partial
 ^{\beta }g_{ij}(y)}\right\vert }\leq \epsilon .$
\end{defin}
We naturally take $\epsilon <1$ in order to have that the riemannian
 metric in the admissible ball be equivalent to the euclidean
 one in ${\mathbb{R}}^{n}.$\ \par 
\quad Of course, without any extra hypotheses on the riemannian manifold
 $M,$ we have $\forall m\in {\mathbb{N}},\ m\geq 2,\ \forall
 \epsilon >0,\ \forall x\in M,$ taking $\displaystyle g_{ij}(x)=\delta
 _{ij}$ in a chart on $\displaystyle B(x,R)$ and the radius $R$
 small enough, the ball $\displaystyle B(x,R)$ is $m,\epsilon
 $ admissible.\ \par 

\begin{defin}
~\label{CL26}Let $\displaystyle x\in M,$ we set $\displaystyle
 R'(x)=\sup \ \lbrace R>0::B(x,R)\ is\ \epsilon \ admissible\rbrace
 .$ We shall say that $\displaystyle R_{\epsilon }(x):=\min \
 (1,R'(x))$ is the $m,\epsilon $ {\bf admissible radius} at $\displaystyle x.$
\end{defin}
\quad Our admissible radius is bigger than the harmonic radius $\displaystyle
 r_{H}(1+\epsilon ,\ m-1,\ 0)$ defined in the Hebey's book ~\cite[p.
 4]{Hebey96}, because we do not require the coordinates to be
 harmonic. I was strongly inspired by this book.\ \par 
\ \par 
\quad When comparing non compact $M$ to the compact case treated above,
 we have four important issues:\ \par 
\quad (0) we have no longer, in general, a global solution $u\in L_{G}^{2}(M)$
 of $Du=\omega $ for a $G$-form $\displaystyle \omega \in L_{G}^{2}(M)$
 verifying $\displaystyle \omega \perp \mathrm{k}\mathrm{e}\mathrm{r}D^{*}.$
 So we have to make this "threshold" hypothesis, which depends on $G.$\ \par 
In case the elliptic operator $D$ is essentially self adjoint,
 this amounts to ask that its spectrum has a gap near $0$: i.e.
 $\exists \delta >0$ such that $D$ has no spectrum in $\rbrack
 0,\delta \lbrack .$ We shall note this hypothesis (THL2G). Moreover,
 because $\displaystyle L^{2}_{G}(M)$ is a Hilbert space, we
 have that the $u\in L^{2}_{G}(M),\ Du=\omega $ with the smallest
 norm is given linearly with respect to $\omega .$ This means
 that the hypothesis (THL2G) gives a bounded linear operator
 $\displaystyle S:\ L^{2}_{G}(M)\rightarrow L^{2}_{G}(M)$ such
 that $D(S\omega )=\omega $ provided that $\displaystyle \omega
 \perp \mathrm{k}\mathrm{e}\mathrm{r}D^{*}.$\ \par 
\ \par 
\quad (1) The "ellipticity constant" may go to zero at infinity and
 we prevent this by asking that $D$ is uniformly elliptic in
 the sense of Definition~\ref{LIR39}.\ \par 
\quad To be sure that the constants in the local elliptic inequalities
 are uniform, we make also the hypothesis that the coefficients
 of $D$ are in ${\mathcal{C}}^{1}(M).$ These are the hypotheses
 (UEAB) in Definition~\ref{mLIR26}.\ \par 
\ \par 
\quad (i) The "admissible" radius may go to $\displaystyle 0$ at infinity,
 which is the case, for instance, if the canonical volume measure
 $\displaystyle dv_{g}$ of $\displaystyle (M,g)$ is finite and
 $M$ is not compact.\ \par 
\quad (ii) If $\displaystyle dv_{g}$ is not finite, which is the case,
 for instance, if the "admissible"  radius is bounded below,
 then $G$-forms in $\displaystyle L^{t}_{G}(M)$ are generally
 not in $\displaystyle L^{r}_{G}(M)$ for $\displaystyle r<t.$\ \par 
\ \par 
\quad We address these two last problems by use of adapted weights
 on $\displaystyle (M,g).$ These weights are relative to a Vitali
 type covering ${\mathcal{C}}_{\epsilon }$ of "admissible balls":
 the weights are positive functions which vary slowly on the
 balls of the covering ${\mathcal{C}}_{\epsilon }.$\ \par 
\quad To state our result in the case of a complete non compact riemannian
 manifold $M$ without boundary we shall use the following definition:\ \par 

\begin{defin}
~\label{HC33}We shall define the Sobolev exponents $\displaystyle
 S_{k}(r)$ by $\displaystyle \ \frac{1}{S_{k}(r)}:=\frac{1}{r}-\frac{k}{n}$
 where $n$ is the dimension of the manifold $M.$
\end{defin}
Now we suppose we have an elliptic operator $D$ with ${\mathcal{C}}^{1}(M)$
 smooth coefficients, of order $m,$ operating on the vector bundle
 $\displaystyle G:=(H,\pi ,M)$ over $M.$  We set $t_{l}:=S_{ml}(2).$
 We suppose that $\displaystyle t_{l-1}\leq r<t_{l},$ and $\displaystyle
 t_{l-1}<\infty .$\ \par 
\quad We set the weights, with $R(x)$ the admissible radius at the
 point $x\in M:$\ \par 
\quad \quad \quad $\displaystyle \ w_{l}(x)=R(x)^{lmt_{l-1}}$ and $\displaystyle
 v_{r}(x):=R(x)^{(\frac{r}{t_{l}}-1)+(l+2)mr}.$\ \par 
\quad Now we can state the main result of this section, where we omit
 the subscript $G$ to ease the notation.\ \par 

\begin{thm}
Under hypotheses (THL2G) and (UEAB), we have, provided that:\par 
\quad \quad \quad $\displaystyle \omega \in L^{2}(M)\cap L^{t_{l-1}}(M,w_{l}),\
 \omega \perp \mathrm{k}\mathrm{e}\mathrm{r}D^{*},$\par 
that $\displaystyle u:=S\omega $ verifies $\displaystyle Du=\omega
 $ with the estimates:\par 
\quad \quad \quad $\displaystyle {\left\Vert{u}\right\Vert}_{L^{r}(M,v_{r})}\leq
 \max ({\left\Vert{\omega }\right\Vert}_{L^{t_{l-1}}(M,w_{l})},{\left\Vert{\omega
 }\right\Vert}_{L^{2}(M)}).$\par 
We also have with the same $u$:\par 
\quad \quad \quad $\displaystyle {\left\Vert{u}\right\Vert}_{W^{m,r}(M,v_{r})}\leq
 c_{1}{\left\Vert{\omega }\right\Vert}_{L^{t_{l}}(M,v_{r})}+c_{2}\max
 ({\left\Vert{\omega }\right\Vert}_{L^{t_{l-1}}(M,w_{l})},{\left\Vert{\omega
 }\right\Vert}_{L^{2}(M)}).$
\end{thm}

\begin{rem}
~\label{LIR32}If the admissible radius $R(x)$ is uniformly bounded
 below, we can forget the weights and we get the existence of
 a solution $u$ of $Du=\omega $ with:\par 
\quad \quad \quad $\displaystyle {\left\Vert{u}\right\Vert}_{L^{r}(M)}\leq \max
 ({\left\Vert{\omega }\right\Vert}_{L^{t_{l-1}}(M)},{\left\Vert{\omega
 }\right\Vert}_{L^{2}(M)}).$\par 
\quad \quad \quad $\displaystyle {\left\Vert{u}\right\Vert}_{W^{m,r}(M)}\leq c_{1}{\left\Vert{\omega
 }\right\Vert}_{L^{t_{l}}(M)}+c_{2}\max ({\left\Vert{\omega }\right\Vert}_{L^{t_{l-1}}(M)},{\left\Vert{\omega
 }\right\Vert}_{L^{2}(M)}).$
\end{rem}
\ \par 
\quad An advantage of this method is that it separates cleanly the
 geometry and the analysis:\ \par 
\quad $\bullet $ The geometry controls the behavior of the admissible
 radius $R(x)$ as a function of $x$ in $M.$ For instance by Theorem
 1.3 in Hebey~\cite{Hebey96}, we have that the harmonic radius
 $\displaystyle r_{H}(1+\epsilon ,\ m,0)$ is bounded below if
 the Ricci curvature $\displaystyle Rc$ verifies $\displaystyle
 \ \forall j\leq m,\ {\left\Vert{\nabla ^{j}Rc}\right\Vert}_{\infty
 }<\infty $  and the injectivity radius is bounded below. This
 implies that the $m,\epsilon $ admissible radius $R(x)$ is also
 bounded below.\ \par 
\quad $\bullet $ The analysis gives the weights as function of $R(x)$
 to get the right estimates. For instance if the admissible radius
 $R(x)$ is bounded below, then we can forget the weights and
 we get more "classical" estimates, as in Remark~\ref{LIR32}.\ \par 
\ \par 
\quad I am indebted to A. Bachelot, B. Helffer, G. M\'etivier and J.
 Sj\"ostrand for clearing strongly my knowledge on the local
 existence of solutions to \emph{system} of elliptic equations
 needed in the study of elliptic equations acting on vector bundles.\ \par 
\quad This work is presented the following way.\ \par 
\quad $\bullet $ In the next section we state the LIR method in the
 general context of metric spaces.\ \par 
\quad $\bullet $ In Section~\ref{LIR33} we apply it for the case of
 elliptic equations in a compact connected riemannian manifold
 without boundary.\ \par 
\quad $\bullet $ In Section~\ref{LIR34} we study the case of elliptic
 equations in a compact connected riemannian manifold with a
 smooth boundary.\ \par 
\quad $\bullet $ In Section~\ref{LIR35} we show that the LIR condition,
 which is a priori estimates, implies the existence of a local
 solution with good estimates.\ \par 
\quad $\bullet $ In section~\ref{LIR36} we study the more delicate
 case of elliptic equations in a complete non compact connected
 riemannian manifold without boundary.\ \par 
\quad $\bullet $ Finally in the Appendix we have put technical results
 concerning the $\epsilon $ admissible balls,  Vitali coverings
 and Sobolev spaces.\ \par 
\ \par 
\quad If the general ideas under this work are quite simple and natural,
 unfortunately the computations to make them work are a little
 bit technical.\ \par 
\ \par 

\section{The Local Increasing Regularity Method (LIRM).}
\quad Let $X$ be a complete metric space with a positive $\sigma $-finite
 measure $\mu .$ Let $\Omega $ be a relatively compact domain
 in $X.$ We shall denote $E^{p}(\Omega )$ the set of ${\mathbb{C}}^{p}$
 valued fonctions on $\Omega .$\ \par 
This means that $\omega \in E^{p}(X)\iff \omega (x)=(\omega _{1}(x),...,\omega
 _{p}(x)).$ We put a punctual norm on $\omega $ in $\displaystyle
 E^{p}(\Omega )$ the following way: for any $x\in \Omega ,\ \left\vert{\omega
 (x)}\right\vert ^{2}:=\sum_{j=1}^{p}{\left\vert{\omega _{j}(x)}\right\vert
 ^{2}}.$ We consider the Lebesgue space $L^{r}_{p}(\Omega ),$ i.e.\ \par 
\quad \quad \quad $\omega \in L^{r}_{p}(\Omega )\iff {\left\Vert{\omega }\right\Vert}^{r}_{L^{r}_{p}(\Omega
 )}:=\int_{\Omega }{\left\vert{\omega (x)}\right\vert ^{r}d\mu
 (x)}<\infty .$\ \par 
The space $L^{2}_{p}(\Omega )$ is a Hilbert space with the scalar
 product ${\left\langle{\omega ,\omega '}\right\rangle}:=\int_{\Omega
 }{{\left({\sum_{j=1}^{p}{\omega _{j}(x)\bar \omega '_{j}(x)}}\right)}d\mu
 (x)}.$\ \par 
\quad We are interested in solutions of a linear equation $Du=\omega
 ,$ where $D=D_{p}$ is a linear operator acting on $E^{p}.$ This
 means that $D$ is a matrix whose entries are linear operators
 on functions.\ \par 
\quad We shall make the following hypotheses.\ \par 
Let $\Omega $ be a relatively compact connected domain in $X.$
 	Let $B:=B(x,R)$ be a ball in $X$ and $B^{1}:=B(x,R/2).$ There
 is a $\tau >0$ with $\displaystyle \frac{1}{t}=\frac{1}{r}-\tau
 $ such that:\ \par 
\quad \quad \quad (i) Local Increasing Regularity (LIR), we have:\ \par 
\ \par 
\quad \quad \quad $\forall x\in \bar \Omega ,\ \exists R>0::\forall r\geq s,\ \exists
 c_{l}>0,\ \forall u\in L_{p}^{r}(B),\ {\left\Vert{u}\right\Vert}_{L_{p}^{t}(B^{1})}\leq
 c_{l}({\left\Vert{Du}\right\Vert}_{L_{p}^{r}(B)}+{\left\Vert{u}\right\Vert}_{L_{p}^{r}(B)}).$\
 \par 
\ \par 
\quad It may happen, in the case $X$ is a manifold, that we have a
 better regularity locally:\ \par 
\quad \quad \quad (i') Local Increasing Regularity (LIR) with Sobolev estimates:
 there is $\alpha >0$ such that\ \par 
\ \par 
\quad \quad \quad $\forall x\in \bar \Omega ,\ \exists R>0::\forall r\geq s,\ \exists
 c_{l}>0,\ \forall u\in L_{p}^{r}(B),\ {\left\Vert{u}\right\Vert}_{W_{p}^{\alpha
 ,r}(B^{1})}\leq c_{l}({\left\Vert{Du}\right\Vert}_{L_{p}^{r}(B)}+{\left\Vert{u}\right\Vert}_{L_{p}^{r}(B)}).$\
 \par 
\ \par 
\quad \quad (ii) Global resolvability. There exists a threshold $s\in (1,\infty
 )$ such that we can solve $Dw=\omega $ globally in $\Omega $
 with $L^{s}-L^{s}$ estimates. It may happen that there is a
 constrain:  let $K$ be a subspace of $L_{p}^{s'}(\Omega ),\
 s'$ the conjugate exponent of $s,$ then we can solve $Dw=\omega
 $ if $\displaystyle \omega \perp K.$ In case with no constrain,
 we set $K=\lbrace 0\rbrace .$ This means:\ \par 
\ \par 
\quad \quad \quad \quad 	$\exists c_{g}>0,\ \exists w\ s.t.\ Dw=\omega $ in $\displaystyle
 \Omega $ and $\ {\left\Vert{w}\right\Vert}_{L_{p}^{s}(\Omega
 )}\leq c_{g}{\left\Vert{\omega }\right\Vert}_{L_{p}^{s}(\Omega
 )},$ provided that $\omega \perp K.$\ \par 
\ \par 
\quad It may happen, in the case $X$ is a manifold, that we have a
 better regularity for the global existence:\ \par 
\quad \quad (ii') Sobolev regularity: We can solve $Dw=\omega $ globally
 in $\Omega $ with $L^{s}-W^{\alpha ,s}$ estimates, i.e.\ \par 
\ \par 
\quad \quad \quad \quad 	$\exists c_{g}>0,\ \exists w\ s.t.\ Dw=\omega $ in $\displaystyle
 \Omega $ and $\ {\left\Vert{w}\right\Vert}_{W_{p}^{\alpha ,s}(\Omega
 )}\leq c_{g}{\left\Vert{\omega }\right\Vert}_{L_{p}^{s}(\Omega
 )},$ provided that $\displaystyle \omega \perp K.$\ \par 
\ \par 
Then we have:\ \par 

\begin{thm}
~\label{LIR1}Under the assumptions (i), (ii) above, there is
 a positive constant $\displaystyle c_{f}$ such that for $r\geq
 s,$ if $\omega \in L_{p}^{r}(\Omega ),\ \omega \perp K$ there
 is a$\ u\in L_{p}^{t}(\Omega )$ with $\displaystyle \frac{1}{t}=\frac{1}{r}-\tau
 ,$ such that $Du=\omega $ and ${\left\Vert{u}\right\Vert}_{L_{p}^{t}(\Omega
 )}\leq c_{f}{\left\Vert{\omega }\right\Vert}_{L_{p}^{r}(\Omega )}.$\par 
\quad If moreover we have (i') and (ii') and the manifold $X$ admits
 the Sobolev embedding theorems, then $u\in W_{p}^{\alpha ,r}(\Omega
 )$ with control of the norm.
\end{thm}
\quad Proof.\ \par 
Let $\omega \in L^{r}_{p}(\Omega ),\ r>s.$ Because $\Omega $
 is relatively compact and $\mu $ is $\sigma $-finite, we have
 that $\displaystyle \omega \in L^{s}_{p}(\Omega ).$ The global
 resolvability, condition (ii), gives that there is a $u\in L^{s}_{p}(\Omega
 )$ such that $Du=\omega ,$ provided that $\displaystyle \omega \perp K.$\ \par 
The LIR, condition (i), gives that, for any $x\in \bar \Omega
 $ there is a ball $B:=B(x,R)$ and a smaller ball $\displaystyle
 B^{1}:=B(x,R/2)$ such that, with $\displaystyle \frac{1}{t_{1}}=\frac{1}{s}-\tau
 $ (we often forget the subscript $p$ for simplicity):\ \par 
\ \par 
\quad \quad \quad $\displaystyle {\left\Vert{u}\right\Vert}_{L^{t_{1}}(B^{1})}\leq
 C({\left\Vert{Du}\right\Vert}_{L^{s}(B)}+{\left\Vert{u}\right\Vert}_{L^{s}(B)})=C({\left\Vert{\omega
 }\right\Vert}_{L^{s}(B)}+{\left\Vert{u}\right\Vert}_{L^{s}(B)})\leq
 C({\left\Vert{\omega }\right\Vert}_{L^{r}(B)}+{\left\Vert{u}\right\Vert}_{L^{s}(B)}),$\
 \par 
\ \par 
because ${\left\Vert{\omega }\right\Vert}_{L^{s}(B)}\lesssim
 {\left\Vert{\omega }\right\Vert}_{L^{r}(B)},$ since $r\geq s$
 and $\bar \Omega $ is compact.\ \par 
Then applying again the LIR we get, with the smaller ball $B^{2}:=B(x,R/4)$
 and with $t_{2}:=\min (r,t_{1}),$\ \par 
\ \par 
\quad \quad \quad ${\left\Vert{u}\right\Vert}_{L^{t_{2}}(B^{2})}\leq C({\left\Vert{\omega
 }\right\Vert}_{L^{t_{1}}(B^{1})}+{\left\Vert{u}\right\Vert}_{L^{t_{1}}(B^{1})})\lesssim
 ({\left\Vert{\omega }\right\Vert}_{L^{r}(B)}+{\left\Vert{u}\right\Vert}_{L^{s}(B)}).$\
 \par 
\ \par 
\quad $\bullet $ If $t_{1}\geq r\Rightarrow t_{2}=r,$ and ${\left\Vert{u}\right\Vert}_{L^{r}(B^{1})}\lesssim
 ({\left\Vert{Du}\right\Vert}_{L^{r}(B)}+{\left\Vert{u}\right\Vert}_{L^{s}(B)})$
 and with $\displaystyle \frac{1}{t}=\frac{1}{r}-\tau ,$\ \par 
\ \par 
\quad \quad \quad ${\left\Vert{u}\right\Vert}_{L^{t}(B^{2})}\lesssim ({\left\Vert{\omega
 }\right\Vert}_{L^{r}(B^{1})}+{\left\Vert{u}\right\Vert}_{L^{r}(B^{1})})\lesssim
 ({\left\Vert{\omega }\right\Vert}_{L^{r}(B)}+{\left\Vert{u}\right\Vert}_{L^{s}(B)}).$\
 \par 
\ \par 
\quad It remains to cover $\bar \Omega $ by a finite set of balls $B^{2}$
 to be done, because\ \par 
\ \par 
\quad \quad \quad $\displaystyle \sum_{B^{2}}{{\left\Vert{u}\right\Vert}_{L^{t}(B)}}\lesssim
 {\left\Vert{u}\right\Vert}_{L^{t}(\Omega )}$ and $\displaystyle
 {\left\Vert{u}\right\Vert}_{L^{s}(\Omega )}\lesssim {\left\Vert{\omega
 }\right\Vert}_{L^{s}(\Omega )}$ by the threshold hypothesis.\ \par 
\ \par 
\quad $\bullet $ If $t_{1}<r,$ we still have:\ \par 
\ \par 
\quad \quad \quad $\displaystyle {\left\Vert{u}\right\Vert}_{L^{t_{2}}(B^{2})}\lesssim
 ({\left\Vert{\omega }\right\Vert}_{L^{r}(B^{1})}+{\left\Vert{u}\right\Vert}_{L^{t_{1}}(B^{1})}).$\
 \par 
\ \par 
Then applying again the LIR we get, with the smaller ball $B^{3}:=B(x,R/8)$
 and with $t_{3}:=\min (r,t_{2}),$\ \par 
\ \par 
\quad \quad \quad $\displaystyle {\left\Vert{u}\right\Vert}_{L^{t_{3}}(B^{3})}\lesssim
 ({\left\Vert{\omega }\right\Vert}_{L^{r}(B^{2})}+{\left\Vert{u}\right\Vert}_{L^{t_{2}}(B^{2})})\lesssim
 ({\left\Vert{\omega }\right\Vert}_{L^{r}(B^{1})}+{\left\Vert{u}\right\Vert}_{L^{t_{1}}(B^{1})})\lesssim
 ({\left\Vert{\omega }\right\Vert}_{L^{r}(B)}+{\left\Vert{u}\right\Vert}_{L^{s}(B)}).$\
 \par 
\ \par 
Hence if $t_{2}\geq r$ we are done as above, if not we repeat
 the process. Because $\frac{1}{t_{k}}=\frac{1}{s}-k\tau $ after
 a finite number $k\leq 1+\frac{1}{\tau }(\frac{r-s}{2s})$ of
 steps we have $t_{k}\geq r$ and we get, with $\displaystyle
 B^{k}:=B(x,R/2^{k})$ and another constant $C,\ {\left\Vert{u}\right\Vert}_{L^{t_{k}}(B^{k})}\leq
 C({\left\Vert{\omega }\right\Vert}_{L^{r}(B)}+{\left\Vert{u}\right\Vert}_{L^{s}(B)}).$\
 \par 
It remains to cover $\bar \Omega $ with a finite number of balls
 $B^{k}(x)$ to prove the first part of the theorem.\ \par 
\ \par 
\quad For the second part, the global resolvability, condition (ii),
 gives that there is a global solution $u\in L^{s}(\Omega )$
 such that $Du=\omega $ in $\Omega $ with ${\left\Vert{u}\right\Vert}_{L^{s}(\Omega
 )}\lesssim {\left\Vert{\omega }\right\Vert}_{L^{s}(\Omega )}.$
 Now if we have the LIR with Sobolev estimates, condition (i'), then\ \par 
\ \par 
\quad \quad \quad $\forall x\in \bar \Omega ,\ \exists R>0::\forall r\geq s,\ \exists
 C>0,\ \forall v\in L^{r}(B(x,R)),\ {\left\Vert{v}\right\Vert}_{W^{\alpha
 ,r}(B^{1})}\leq C({\left\Vert{Dv}\right\Vert}_{L^{r}(B)}+{\left\Vert{v}\right\Vert}_{L^{r}(B)})$\
 \par 
with, as usual, $B:=B(x,R)$ and $B^{1}:=B(x,R/2).$\ \par 
\ \par 
So, because $r\geq s,$ and $\bar \Omega $ is compact, $\displaystyle
 \omega \in L^{s}(\Omega )$ and we get\ \par 
\ \par 
\quad \quad \quad ${\left\Vert{u}\right\Vert}_{W^{\alpha ,s}(B^{1})}\lesssim ({\left\Vert{Du}\right\Vert}_{L^{s}(B)}+{\left\Vert{u}\right\Vert}_{L^{s}(B)})\lesssim
 ({\left\Vert{\omega }\right\Vert}_{L^{r}(B)}+{\left\Vert{u}\right\Vert}_{L^{s}(B)}).$\
 \par 
\ \par 
The Sobolev embedding theorems, true by assumption here, give
 ${\left\Vert{u}\right\Vert}_{L^{\tau }(B^{1})}\leq c{\left\Vert{u}\right\Vert}_{W^{\alpha
 ,s}(B^{1})}$ with $\frac{1}{\tau }=\frac{1}{s}-\frac{\alpha }{n}.$\ \par 
\quad So applying again the LIR condition in a ball $B^{2}:=B(x,R/4),$
 we get, with $t_{1}:=\min (\tau ,r),$\ \par 
\ \par 
\quad \quad \quad $\displaystyle {\left\Vert{u}\right\Vert}_{W^{\alpha ,t_{1}}(B^{2})}\lesssim
 ({\left\Vert{\omega }\right\Vert}_{L^{t_{1}}(B)}+{\left\Vert{u}\right\Vert}_{L^{t_{1}}(B^{1})})\lesssim
 ({\left\Vert{\omega }\right\Vert}_{L^{r}(B)}+{\left\Vert{u}\right\Vert}_{L^{s}(B)}).$\
 \par 
\ \par 
Now we proceed as above. If $\tau \geq r\Rightarrow t_{1}=r,$
 then we apply again the LIR condition to a smaller ball $B^{3}:=B(x,R/8),$
 we get\ \par 
\ \par 
\quad \quad \quad $\displaystyle {\left\Vert{u}\right\Vert}_{W^{\alpha ,r}(B^{3})}\lesssim
 ({\left\Vert{\omega }\right\Vert}_{L^{r}(B)}+{\left\Vert{u}\right\Vert}_{L^{r}(B^{2})})\lesssim
 ({\left\Vert{\omega }\right\Vert}_{L^{r}(B)}+{\left\Vert{u}\right\Vert}_{L^{s}(B)}).$\
 \par 
\ \par 
and we are done by covering $\bar \Omega $ by a finite set of
 balls $B^{3}$ as above.\ \par 
\quad If $\tau <r,$ then we iterate the process as in the previous
 part, adding the use of the Sobolev embedding theorem to increase
 the exponent, up to the moment we reach $r.$ $\hfill\blacksquare $\ \par 

\begin{rem}
~\label{LIR31} We notice that in fact the solution $u$ in Theorem~\ref{LIR1}
 is the same as the one given by condition (ii). It is a case
 of "self improvement" of estimates.
\end{rem}

\section{Application to elliptic PDE.~\label{LIR33}}
\quad Let $(M,g)$ be a ${\mathcal{C}}^{\infty }$ smooth connected compact
 riemannian manifold without boundary. We shall denote $\displaystyle
 G:=(H,\pi ,M)$ a complex ${\mathcal{C}}^{m}$ vector bundle over
 $M$ of rank $N$ with fiber $H.$ The fiber $\pi ^{-1}(x)\simeq
 H$ is equipped with a scalar product varying smoothly with $x$ in $M.$ \ \par 
\ \par 
We can define punctually, for $\omega ,\varphi \in {\mathcal{C}}^{\infty
 }_{G}(M),$ two smooth sections of $G$ over $M,$ a scalar product
 $(\omega ,\varphi )(x):={\left\langle{\omega (x),\varphi (x)}\right\rangle}_{H_{x}}$
 where $H_{x}:=\pi ^{-1}(x)$ is the fiber over $x\in M.$ This
 gives a modulus: for $x\in M,\ \left\vert{\omega }\right\vert
 (x):={\sqrt{(\omega ,\omega )(x)}}.$ By use of the canonical
 volume $dv_{g}$ on $M$ we get a scalar product:\ \par 
\quad \quad \quad $\displaystyle {\left\langle{\omega ,\varphi }\right\rangle}:=\int_{M}{(\omega
 ,\varphi )(x)dv_{g}(x)},$\ \par 
for $G$-forms in $L^{2}_{G}(M)$ i.e. such that\ \par 
\quad \quad \quad $\displaystyle {\left\Vert{\omega }\right\Vert}_{L^{2}_{G}(M)}^{2}:=\int_{M}{\left\vert{\omega
 }\right\vert ^{2}(x)dv_{g}(x)}<\infty .$\ \par 
The same way we define the spaces $L^{r}_{G}(M)$ of $G$-forms
 $\omega $ such that:\ \par 
\quad \quad \quad $\displaystyle {\left\Vert{\omega }\right\Vert}_{L^{r}_{G}(M)}^{2}:=\int_{M}{\left\vert{\omega
 }\right\vert ^{2}(x)dv_{g}(x)}<\infty .$\ \par 
\quad Let $D:G\rightarrow G$ be a linear differential operator of order
 $m$ with ${\mathcal{C}}^{\infty }$ coefficients. There is a
 formal adjoint $\displaystyle D^{*}:G\rightarrow G$ defined
 by the identity ${\left\langle{D^{*}f,g}\right\rangle}={\left\langle{f,Dg}\right\rangle}.$\
 \par 
\ \par 
\quad We shall use the definition of ellipticity given by Warner~\cite[Definition
 6.28, p. 240]{Warner83} or by Donaldson~\cite[ p. 17]{Donaldson08}.\ \par 
\quad Let $D:E\rightarrow F$ be a differential operator of order $m$
 operating from the sections of the vector bundle $E$ to the
 ones of the vector bundle $F$ over $M.$ Then at each point $x\in
 M$ and for each cotangent vector $\xi \in T^{*}M$ there is a
 linear map $\sigma _{\xi }:E_{x}\rightarrow F_{x}$ which can
 be defined the following way:\ \par 
choose a section $s$ of $E,$ and a function $f$ on $M,$ vanishing
 at $x$ and with $df=\xi $ at $x.$ Then we can define       
     $\sigma _{\xi }(s(x))=D(f^{m}s)(x).$ We can check that this
 definition is independent of the choice of $f,s.$ Now we can state:\ \par 

\begin{defin}
~\label{LIR39} An operator $D:E\rightarrow F$ is elliptic if
 for each nonzero $\xi \in TM_{x},$ the linear map $\sigma _{\xi
 }$ is an isomorphism from $E_{x}$ to $\displaystyle F_{x}.$
 We shall say that $D$ is uniformly elliptic if the isomorphism
 $\sigma _{\xi }$ and its inverse are bounded independently of
 the point $x\in M$ for $\left\vert{\xi }\right\vert =1.$
\end{defin}
Then for $s=2,$ Warner~\cite[Exercice 21, p. 257]{Warner83} or
 also Donaldson~\cite[Theorem 4, p. 16]{Donaldson08}, proved:\ \par 

\begin{thm}
~\label{LIR4} Let $D$ be an operator of order $m$ acting on sections
 of $G:=(H,\pi ,M)$ in the connected compact riemannian manifold
 $M$ without boundary. Suppose that $D$ is elliptic and with
 ${\mathcal{C}}^{\infty }$ smooth coefficients. \par 
\quad 1. In $L^{2}_{G}(M),\ \mathrm{k}\mathrm{e}\mathrm{r}D,\ \mathrm{k}\mathrm{e}\mathrm{r}D^{*}$
 are finite dimensional vector spaces.\par 
\quad 2. We can solve the equation $Du=\omega $ in $\displaystyle L^{2}_{G}(M)$
 if and only if $\omega $ is orthogonal to $\mathrm{k}\mathrm{e}\mathrm{r}D^{*}.$
\end{thm}

      Moreover, because $\displaystyle L^{2}_{G}(M)$ is a Hilbert
 space, we have that there is a bounded linear operator $\displaystyle
 S:\ L^{2}_{G}(M)\rightarrow L^{2}_{G}(M)$ such that $D(S\omega
 )=\omega $ provided that $\displaystyle \omega \perp \mathrm{k}\mathrm{e}\mathrm{r}D^{*}.$\
 \par 
\ \par 
\quad On the other hand we have local interior regularity by H\"ormander~\cite[Theorem
 17.1.3, p. 6]{Hormand94}, in the case of functions. We quote
 it in the weakened form we need:\ \par 

\begin{thm}
(LIR) Let $D$ be an operator of order $m$ on ${\mathcal{C}}^{\infty
 }(M)$ in the complete riemannian manifold $M.$ Suppose that
 $D$ is elliptic and with ${\mathcal{C}}^{\infty }$ smooth coefficients.
 Then, for any $x\in M$ there is a ball $B_{x}:=B(x,R)$ and a
 smaller ball $\displaystyle B'_{x}$ relatively compact in $B_{x},$
 such that:\par 
\quad \quad \quad ${\left\Vert{u}\right\Vert}_{W^{m,r}(B'_{x})}\leq C({\left\Vert{Du}\right\Vert}_{L^{r}(B_{x})}+{\left\Vert{u}\right\Vert}_{L^{r}(B_{x})}).$
\end{thm}

      For the case of $G$-forms, we need to use Agmon, Douglis
 and Nirenberg~\cite[Theorem 10.3]{AgDougNir64}:\ \par 

\begin{thm}
~\label{LIR2}Positive constants $r_{1}$ and $K_{1}$ exist such
 that, if $r\leq r_{1}$ and the ${\left\Vert{u_{j}}\right\Vert}_{t_{j}},\
 j=1,...,N,$ are finite, then $\displaystyle {\left\Vert{u_{j}}\right\Vert}_{l+t_{j}}$
 also is finite for $\displaystyle j=1,...,N,$ and\par 
\quad \quad \quad $\displaystyle {\left\Vert{u_{j}}\right\Vert}_{l+t_{j}}\leq K_{1}{\left({\sum_{j}{{\left\Vert{F_{j}}\right\Vert}_{l-s_{j}}}+\sum_{j}{{\left\Vert{u_{j}}\right\Vert}_{0}}}\right)}.$\par
 
The constants $r_{1},K_{1}$ depend on $n,N,t',A,b,p,k,$ and $l$
 and also on the modulus of continuity of the leading coefficients
 in the $l_{ij}.$
\end{thm}
\quad From this theorem we get quite easily what we want (in the case
 $r=2$ and in its \emph{global version}, F.W. Warner~\cite[Theorem
 6.29, p. 240]{Warner83} quotes it as \emph{Fundamental Inequality}):\ \par 

\begin{thm}
~\label{LIR0}(LIR) Let $D$ be an operator of order $m$ on $G$
 in the complete riemannian manifold $M.$ Suppose that $D$ is
 elliptic and with ${\mathcal{C}}^{1}(M)$ smooth coefficients.
 Then, for any $x\in M$ there is a ball $B:=B(x,R)$ and, with
 the ball $\displaystyle B^{1}:=B(x,R/2),$ we have:\par 
\quad \quad \quad ${\left\Vert{u}\right\Vert}_{W^{m,r}_{G}(B^{1})}\leq c_{1}{\left\Vert{Du}\right\Vert}_{L^{r}_{G}(B)}+c_{2}R^{-m}{\left\Vert{u}\right\Vert}_{L^{r}_{G}(B)}).$\par
 
Moreover the constants are independent of the radius $R$ of the ball $B.$
\end{thm}
\quad Proof.\ \par 
Let $\displaystyle x\in M$ we choose a chart $(V,\ \varphi (y))$
 so that $g_{ij}(x)=\delta _{ij}$ and $\varphi (V)=B_{e}$ where
 $B_{e}=B_{e}(0,R_{e})$ is a Euclidean ball centered at $\varphi
 (x)=0$ and $g_{ij}$ are the components of the metric tensor
 w.r.t. $\varphi .$ We choose also the chart $\displaystyle (V,\
 \varphi )$ to trivialise the bundle $G.$ So read in $\displaystyle
 (V,\varphi )$ we have that the sections of $G$ are just ${\mathbb{C}}^{N}$
 valued functions.\ \par 
\quad We denote by $D_{\varphi }$ the operator $D$ read in the map
 $(V,\varphi ).$ This is still an elliptic system operating on
 ${\mathbb{C}}^{N}$ valued functions in $B_{e}$ in ${\mathbb{R}}^{n}.$
 Let $\chi \in {\mathcal{D}}(B_{e})$ such that $\chi =1$ in $B^{1}_{e}:=B_{e}(0,R_{e}/2)\Subset
 B_{e}.$ Let $u$ be a $G$-form in $L^{r}_{G}(\varphi ^{-1}(B_{e}))$
 such that $Du$ is also in $\displaystyle L^{r}_{G}(\varphi ^{-1}(B_{e})).$
 Denote by $u_{\varphi }$ the ${\mathbb{C}}^{N}$ valued functions
 $u$ read in $(V,\varphi ).$ We can apply the Agmon, Douglis
 and Nirenberg Theorem~\ref{LIR2} to $\chi u_{\varphi }$ and
 we get, with the constant $K$ independent of the radius $R_{e}$
 of $B_{e},$\ \par 
\quad \quad \quad \begin{equation} {\left\Vert{ \chi u_{\varphi }}\right\Vert}_{W^{m,r}(B_{e})}\leq
 K({\left\Vert{D_{\varphi }(\chi u_{\varphi })}\right\Vert}_{L^{r}(B_{e})}+R_{e}^{-m}{\left\Vert{\chi
 u_{\varphi }}\right\Vert}_{L^{r}(B_{e})}).\label{LIR3}\end{equation}\ \par 
We have that $D_{\varphi }(\chi u_{\varphi })=\chi D_{\varphi
 }(u_{\varphi })+u_{\varphi }D_{\varphi }\chi +\Delta _{\varphi
 },$ with $\Delta _{\varphi }:=D_{\varphi }(\chi u_{\varphi })-\chi
 D_{\varphi }(u_{\varphi })-u_{\varphi }D_{\varphi }\chi .$ The
 point is that $\Delta _{\varphi }$ contains only derivatives
 of the $j^{th}$ component of $u_{\varphi }$ of order strictly
 less than in the $j^{th}$ component of $u_{\varphi }$ in $D_{\varphi
 }u_{\varphi }.$ So we have\ \par 
\ \par 
\quad \quad \quad ${\left\Vert{\Delta _{\varphi }}\right\Vert}_{L^{r}(B_{e})}\leq
 {\left\Vert{\partial \chi }\right\Vert}_{\infty }{\left\Vert{\chi
 u_{\varphi }}\right\Vert}_{W^{m-1,r}(B_{e})}\leq R_{e}^{-1}{\left\Vert{\chi
 u_{\varphi }}\right\Vert}_{W^{m-1,r}(B_{e})}.$\ \par 
\ \par 
We can use the "Peter-Paul" inequality~\cite[Theorem 7.28, p.
 173]{GuilbargTrudinger98} (see also~\cite[Theorem 6.18, (g)
 p. 232]{Warner83} for the case $r=2.$)\ \par 
\quad \quad \quad $\displaystyle \forall \epsilon >0,\ \exists C_{\epsilon }>0::{\left\Vert{\chi
 u_{\varphi }}\right\Vert}_{W^{m-1,r}(B_{e})}\leq \epsilon {\left\Vert{\chi
 u_{\varphi }}\right\Vert}_{W^{m,r}(B_{e})}+C\epsilon ^{-m+1}{\left\Vert{\chi
 u_{\varphi }}\right\Vert}_{L^{r}(B_{e})}.$\ \par 
We choose $\epsilon =R_{e}\eta $ and we get\ \par 
\quad \quad \quad $\displaystyle R_{e}^{-1}{\left\Vert{\chi u_{\varphi }}\right\Vert}_{W^{m-1,r}(B_{e})}\leq
 \eta {\left\Vert{\chi u_{\varphi }}\right\Vert}_{W^{m,r}(B_{e})}+C\eta
 ^{-m+1}R_{e}^{-m}{\left\Vert{\chi u_{\varphi }}\right\Vert}_{L^{r}(B_{e})}.$\
 \par 
Putting this in~(\ref{LIR3}) we get\ \par 
\ \par 
\quad $\displaystyle {\left\Vert{\chi u_{\varphi }}\right\Vert}_{W^{m,r}(B_{e})}\leq
 K({\left\Vert{\chi D_{\varphi }u_{\varphi }}\right\Vert}_{L^{r}(B_{e})}+\eta
 {\left\Vert{\chi u_{\varphi }}\right\Vert}_{W^{m,r}(B_{e})}+$\ \par 
\quad \quad \quad \quad \quad \quad \quad \quad \quad \quad \quad \quad \quad $\displaystyle +C\eta ^{-m+1}R_{e}^{-m}{\left\Vert{\chi u_{\varphi
 }}\right\Vert}_{L^{r}(B_{e})}+{\left\Vert{u_{\varphi }D_{\varphi
 }\chi }\right\Vert}_{L^{r}(B_{e})}).$\ \par 
\ \par 
But again ${\left\Vert{D_{\varphi }\chi }\right\Vert}_{\infty
 }\leq R_{e}^{-m}$ so, choosing $\eta $ small enough to get $\eta
 K\leq 1/2,$ we have with new constants still independent of $R_{e}$:\ \par 
\quad \quad \quad $\frac{1}{2}{\left\Vert{\chi u_{\varphi }}\right\Vert}_{W^{m,r}(B_{e})}\leq
 c_{1}{\left\Vert{\chi D_{\varphi }u_{\varphi }}\right\Vert}_{L^{r}(B_{e})}+c_{2}R_{e}^{-m}{\left\Vert{\chi
 u_{\varphi }}\right\Vert}_{L^{r}(B_{e})}.$\ \par 
\ \par 
Now $\chi =1$ in $B^{1}_{e}$ and $\chi \leq 1$ gives, changing
 the constants suitably:\ \par 
\quad \quad \quad \begin{equation} {\left\Vert{ u_{\varphi }}\right\Vert}_{W^{m,r}(B^{1}_{e})}\leq
 c_{1}{\left\Vert{D_{\varphi }u_{\varphi }}\right\Vert}_{L^{r}(B_{e})}+c_{2}R_{e}^{-m}{\left\Vert{u_{\varphi
 }}\right\Vert}_{L^{r}(B_{e})}.\label{LIR37}\end{equation}\ \par 
\quad It remains to go back to the manifold $M$ to end the proof. $\hfill\blacksquare
 $\ \par 
\ \par 
\quad We deduce the local\emph{ elliptic inequalities:}\ \par 

\begin{cor}
~\label{LIR6} Let $D$ be an operator of order $m$ on $G$ in the
 complete riemannian manifold $M.$ Suppose that $D$ is elliptic
 and with ${\mathcal{C}}^{1}(M)$ smooth coefficients. Then, for
 any $x\in M$ there is a ball $B:=B(x,R)$ and the smaller ball
 $\displaystyle B^{1}:=B(x,R/2),$ such that, $\forall k\in {\mathbb{N}},$
 with $D$ in ${\mathcal{C}}^{k+1}(M)$ here, we get for any $G$-form
  $\displaystyle u\in W^{m+k,r}_{G}(B^{1}):$\par 
\quad \quad \quad $\displaystyle {\left\Vert{u}\right\Vert}_{W^{m+k,r}_{G}(B^{1})}\leq
 \sum_{j=0}^{k}{c_{j}R^{-jm}{\left\Vert{Du}\right\Vert}}_{W^{k-j,r}_{G}(B)}+c_{k+1}R^{-(k+1)m}{\left\Vert{u}\right\Vert}_{L^{r}_{G}(B)}.$\par
 
Moreover the constants are independent of the radius $R$ of the ball $B.$
\end{cor}
\quad Proof.\ \par 
As for Theorem~\ref{LIR0}, we choose a chart $(V,\ \varphi )$
 trivialising the bundle $G$ and so that $g_{ij}(x)=\delta _{ij}$
 and $\varphi (V)=B$ where $B$ is a Euclidean ball centered at
 $\varphi (x)=0$ and $g_{ij}$ are the components of the metric
 tensor w.r.t. $\varphi .$ We start with the equation~(\ref{LIR37})
 in ${\mathbb{R}}^{n}$ and we apply it to $\partial _{j}u_{\varphi
 }:=\frac{\partial u_{\varphi }}{\partial y_{j}}$ instead of
 $u_{\varphi }.$ We get\ \par 
\quad \quad \quad $\displaystyle {\left\Vert{\partial _{j}u_{\varphi }}\right\Vert}_{W^{m,r}(B^{1})}\leq
 c_{1}{\left\Vert{D_{\varphi }(\partial _{j}u_{\varphi })}\right\Vert}_{L^{r}(B)}+c_{2}R_{e}^{-m}{\left\Vert{\partial
 _{j}u_{\varphi }}\right\Vert}_{L^{r}(B)}.$\ \par 
\ \par 
Now $\displaystyle D_{\varphi }(\partial _{j}u_{\varphi })=\partial
 _{j}D_{\varphi }(u_{\varphi })+\lbrack D_{\varphi },\partial
 _{j}\rbrack u_{\varphi },$ with as usual, $\displaystyle \lbrack
 D_{\varphi },\partial _{j}\rbrack u_{\varphi }:=D_{\varphi }(\partial
 _{j}u_{\varphi })-\partial _{j}D_{\varphi }(u_{\varphi }).$\ \par 
So we get\ \par 
\ \par 
\quad \quad \quad $\displaystyle {\left\Vert{\partial _{j}u_{\varphi }}\right\Vert}_{W^{m,r}(B^{1})}\leq
 c_{1}{\left\Vert{\partial _{j}D_{\varphi }u_{\varphi }}\right\Vert}_{L^{r}(B)}+c_{1}{\left\Vert{\lbrack
 D_{\varphi },\partial _{j}\rbrack u_{\varphi }}\right\Vert}_{L^{r}(B)}+c_{2}R_{e}^{-m}{\left\Vert{\partial
 _{j}u_{\varphi }}\right\Vert}_{L^{r}(B)}.$\ \par 
\ \par 
So, because $\displaystyle \lbrack D_{\varphi },\partial _{j}\rbrack
 $ is a differential operator of order $m,$ we get\ \par 
\ \par 
\quad \quad \quad $\displaystyle {\left\Vert{\partial _{j}u_{\varphi }}\right\Vert}_{W^{m,r}(B^{1})}\leq
 c_{1}{\left\Vert{D_{\varphi }u_{\varphi }}\right\Vert}_{W^{1,r}(B)}+c_{1}{\left\Vert{u_{\varphi
 }}\right\Vert}_{W^{m,r}(B)}+c_{2}R_{e}^{-m}{\left\Vert{u_{\varphi
 }}\right\Vert}_{W^{1,r}(B)}.$\ \par 
\ \par 
This is true for any $j=1,...,n$ so\ \par 
\ \par 
\quad \quad \quad $\displaystyle {\left\Vert{u_{\varphi }}\right\Vert}_{W^{m+1,r}(B^{1})}\leq
 c_{1}{\left\Vert{D_{\varphi }u_{\varphi }}\right\Vert}_{W^{1,r}(B)}+c_{1}{\left\Vert{u_{\varphi
 }}\right\Vert}_{W^{m,r}(B)}+c_{2}R_{e}^{-m}{\left\Vert{u_{\varphi
 }}\right\Vert}_{W^{1,r}(B)}.$\ \par 
\ \par 
We always have $\ {\left\Vert{u_{\varphi }}\right\Vert}_{W^{1,r}(B)}\leq
 {\left\Vert{u_{\varphi }}\right\Vert}_{W^{m,r}(B)}$ hence, with
 other constants $c_{j},$\ \par 
\ \par 
\quad \quad \quad $\displaystyle {\left\Vert{u_{\varphi }}\right\Vert}_{W^{m+1,r}(B^{1})}\leq
 c_{1}{\left\Vert{D_{\varphi }u_{\varphi }}\right\Vert}_{W^{1,r}(B)}+(c_{1}+c_{2}R^{-m}){\left\Vert{u_{\varphi
 }}\right\Vert}_{W^{m,r}(B)}\leq $\ \par 
\quad \quad \quad \quad \quad \quad \quad \quad \quad \quad \quad $\displaystyle \leq c_{1}{\left\Vert{D_{\varphi }u_{\varphi }}\right\Vert}_{W^{1,r}(B)}+c_{2}R^{-m}{\left\Vert{u_{\varphi
 }}\right\Vert}_{W^{m,r}(B)},$\ \par 
\ \par 
because $R\leq 1.$\ \par 
Now we use again equation~(\ref{LIR37}) to get\ \par 
\ \par 
\quad \quad \quad $\displaystyle {\left\Vert{u_{\varphi }}\right\Vert}_{W^{m,r}(B)}\leq
 c_{1}{\left\Vert{D_{\varphi }u_{\varphi }}\right\Vert}_{L^{r}(B)}+c_{2}R_{e}^{-m}{\left\Vert{u_{\varphi
 }}\right\Vert}_{L^{r}(B)}$\ \par 
\ \par 
hence, still with different constants from line to line\ \par 
\ \par 
\quad \quad \quad $\displaystyle {\left\Vert{u_{\varphi }}\right\Vert}_{W^{m+1,r}(B^{1})}\leq
 c_{1}{\left\Vert{D_{\varphi }u_{\varphi }}\right\Vert}_{W^{1,r}(B)}+c_{2}R_{e}^{-m}({\left\Vert{D_{\varphi
 }u_{\varphi }}\right\Vert}_{L^{r}(B)}+c_{2}R_{e}^{-m}{\left\Vert{u_{\varphi
 }}\right\Vert}_{L^{r}(B)})\leq $\ \par 
\quad \quad \quad \quad \quad \quad \quad \quad \quad \quad \quad $\displaystyle \leq c_{1}{\left\Vert{D_{\varphi }u_{\varphi }}\right\Vert}_{W^{1,r}(B)}+c_{2}R_{e}^{-m}{\left\Vert{D_{\varphi
 }u_{\varphi }}\right\Vert}_{L^{r}(B)}+c_{3}R_{e}^{-2m}{\left\Vert{u_{\varphi
 }}\right\Vert}_{L^{r}(B)}).$\ \par 
\ \par 
Now, proceeding by induction along the same lines, we get\ \par 
\ \par 
\quad \quad \quad $\displaystyle {\left\Vert{u_{\varphi }}\right\Vert}_{W^{m+k,r}(B^{1})}\leq
 \sum_{j=0}^{k}{c_{j}R^{-jm}{\left\Vert{D_{\varphi }u_{\varphi
 }}\right\Vert}}_{W^{k-j,r}(B)}+c_{k+1}R_{e}^{-(k+1)m}{\left\Vert{u_{\varphi
 }}\right\Vert}_{L^{r}(B)}.$\ \par 
\ \par 
\quad It remains to go back to the manifold $M$ to end the proof. $\hfill\blacksquare
 $\ \par 

\begin{rem}
We stress here the dependence in $R$ because we shall need it
 to study the case of non compact riemannian manifolds.
\end{rem}
\quad Now we can prove\ \par 

\begin{thm}
Let $(M,g)$ be a ${\mathcal{C}}^{\infty }$ smooth compact riemannian
 manifold without boundary. Let $D:G\rightarrow G$ be an elliptic
 linear differential operator of order $m$ with ${\mathcal{C}}^{\infty
 }(M)$ coefficients.  Let $\omega \in L^{r}_{G}(M)\cap \mathrm{(}\mathrm{k}\mathrm{e}\mathrm{r}D^{*}\mathrm{)}^{\perp
 }$ with $r\geq 2.$ Then there is a $u\in W^{m,r}_{G}(M)$ such
 that $Du=\omega $ on $M.$ Moreover $u$ is given linearly w.r.t. to $\omega .$
\end{thm}
\quad Proof.\ \par 
Let $\omega \in L^{r}_{G}(M)\cap \mathrm{(}\mathrm{k}\mathrm{e}\mathrm{r}D^{*}\mathrm{)}^{\perp
 }$ with $r\geq 2.$ Because $M$ is compact, we have $\displaystyle
 \omega \in L^{2}_{G}(M).$ Theorem~\ref{LIR4} gives us the Global
 Resolvability, condition (ii), with the threshold $s=2,$ and
 with $K:=\mathrm{k}\mathrm{e}\mathrm{r}D^{*},$ i.e. provided
 that $\omega \perp K$:\ \par 
\ \par 
\quad \quad \quad $u:=S\omega \in L^{2}_{G}(M)::Du=\omega ,\ {\left\Vert{u}\right\Vert}_{2}\leq
 C{\left\Vert{\omega }\right\Vert}_{2}.$\ \par 
\ \par 
The Theorem~\ref{LIR0} of Agmon, Douglis and Nirenberg gives
 us the Local Interior Regularity with the Sobolev estimates
 for $\alpha =m.$\ \par 
\quad So we can apply Theorem~\ref{LIR1} and we use Remark~\ref{LIR31}
 to have that $u=S\omega $ so $u$ is given linearly w.r.t. to
 $\omega .$ The proof is complete. $\hfill\blacksquare $\ \par 
\ \par 
\quad By duality we get the range $r<2.$ We shall proceed as we did
 in~\cite{AmarSt13}, using an avatar of the Serre duality~\cite{Serre55}.\ \par 
\quad Let $g\in L^{r'}_{G}(M)\cap \mathrm{k}\mathrm{e}\mathrm{r}D^{\perp
 },$ because $D^{*}$ has the same elliptic properties than $D,$
 we can solve $D^{*}v=g,$ with $r'<2$ and $r'$ conjugate to $r$
 the following way.\ \par 
We know by the previous part that:\ \par 
\quad \quad \quad \begin{equation}  \forall \omega \in L^{r}_{G}(M)\cap \mathrm{(}\mathrm{k}\mathrm{e}\mathrm{r}D^{*}\mathrm{)}^{\perp
 },\ \exists u\in L^{r}_{G}(M),\ Du=\omega .\label{HD5}\end{equation}\ \par 
Consider the linear form\ \par 
\quad \quad \quad $\forall \omega \in L^{r}_{G}(M),\ {\mathcal{L}}(\omega ):={\left\langle{u,g}\right\rangle},$
 \ \par 
where $u$ is a solution of~(\ref{HD5}); in order for ${\mathcal{L}}(\omega
 )$ to be well defined, we need that if $u'$ is another solution
 of $Du'=\omega ,$ then ${\left\langle{u-u',g}\right\rangle}=0;$
 hence we need that $g$ must be "orthogonal" to $G$-forms $\varphi
 $ such that $D\varphi =0,$ which is precisely our assumption.\ \par 
\quad Hence we have that ${\mathcal{L}}(f)$ is well defined and linear;
 moreover\ \par 
\ \par 
\quad \quad \quad $\displaystyle \ \left\vert{{\mathcal{L}}(f)}\right\vert \leq
 {\left\Vert{u}\right\Vert}_{L^{r}(M)}{\left\Vert{g}\right\Vert}_{L^{r'}(M)}\leq
 c{\left\Vert{\omega }\right\Vert}_{L^{r}(M)}{\left\Vert{g}\right\Vert}_{L^{r'}(M)}.$\
 \par 
\ \par 
So this linear form is continuous on $\omega \in L^{r}_{G}(M)\cap
 (\mathrm{k}\mathrm{e}\mathrm{r}D^{*})^{\perp }.$ By the Hahn
 Banach Theorem there is a form $v\in L^{r'}_{G}(M)$ such that:\ \par 
\ \par 
\quad \quad \quad $\displaystyle \forall \omega \in L^{r}_{G}(M)\cap (\mathrm{k}\mathrm{e}\mathrm{r}D^{*})^{\perp
 },\ {\mathcal{L}}(\omega )={\left\langle{\omega ,v}\right\rangle}={\left\langle{u,g}\right\rangle}.$\
 \par 
\ \par 
But $\omega =Du,$ so we have ${\left\langle{\omega ,v}\right\rangle}={\left\langle{Du,v}\right\rangle}={\left\langle{u,D^{*}v}\right\rangle}={\left\langle{u,g}\right\rangle},$
 for any $u\in {\mathcal{C}}^{\infty }_{G}(M).$ Hence we solved
 $D^{*}v=g$ in the sense of distributions with $v\in L^{r'}_{G}(M).$
 So we proved:\ \par 

\begin{thm}
For any $r,\ 1<r\leq 2,$ if $g\in L_{G}^{r}(M)\cap (\mathrm{k}\mathrm{e}\mathrm{r}D)^{\perp
 }$ there is a $v\in L^{r}_{G}(M)$ such that $D^{*}v=g,\ {\left\Vert{v}\right\Vert}_{L^{r}_{G}(M)}\leq
 c{\left\Vert{g}\right\Vert}_{L^{r}_{G}(M)}.$\par 
Moreover the solution is in $W^{m,r}_{G}(M).$
\end{thm}
\quad It remains to prove the "moreover" and for this we use the LIR
 Theorem~\ref{LIR0}: for any $x\in M$ there is a ball $B:=B(x,R)$
 and, with the ball $\displaystyle B^{1}:=B(x,R/2),$ we get:\ \par 
\ \par 
\quad \quad \quad ${\left\Vert{u}\right\Vert}_{W^{m,r}_{G}(B^{1})}\leq C({\left\Vert{Du}\right\Vert}_{L_{G}^{r}(B)}+{\left\Vert{u}\right\Vert}_{L_{G}^{r}(B)}).$\
 \par 
\ \par 
We cover $M$ with a finite number of balls $B^{1}$ to prove the
 theorem. $\hfill\blacksquare $\ \par 
\ \par 
Set ${\mathcal{H}}_{G}^{2}:=\mathrm{k}\mathrm{e}\mathrm{r}D^{*}\cap
 L^{2}_{G}(M).$\ \par 
\quad Because $D$ and $D^{*}$ have the same elliptic properties, we
 finally proved:\ \par 

\begin{thm}
Let $(M,g)$ be a ${\mathcal{C}}^{\infty }$ smooth compact riemannian
 manifold without boundary. Let $D:G\rightarrow G$ be an elliptic
 linear differential operator of order $m$ with ${\mathcal{C}}^{1}$
 coefficients.  Let $\omega \in L^{r}_{G}(M)\cap ({\mathcal{H}}_{G}^{2})^{\perp
 }$ with $r>1.$ Then there is a $u\in L^{r}_{G}(M)$ such that
 $Du=\omega $ on $M.$ Moreover the solution is in $W^{m,r}_{G}(M).$
\end{thm}
\quad Now we make the hypothesis that $D$ has ${\mathcal{C}}^{\infty
 }$ smooth coefficients. The Theorem~\ref{LIR4} of Warner or
 Donaldson gives, on a compact manifold $M$ without boundary,
 that $\mathrm{d}\mathrm{i}\mathrm{m}_{{\mathbb{R}}}{\mathcal{H}}_{G}^{2}<\infty
 .$\ \par 
We shall generalise here a well known result valid for the Hodge
 laplacian.\ \par 

\begin{lem}
~\label{LIR7}We have $\displaystyle {\mathcal{H}}_{G}^{2}\subset
 {\mathcal{C}}^{\infty }(M).$
\end{lem}
\quad Proof.\ \par 
Take $x\in M,\ h\in {\mathcal{H}}_{G}^{2}.$ The fundamental inequalities,
 Corollary~\ref{LIR6}, gives, applied to $D^{*},$ that there
 is a ball $B:=B(x,R)$ with the ball $\displaystyle B^{1}:=B(x,R/2)$
 such that:\ \par 
\ \par 
\quad \quad \quad $\forall k\in {\mathbb{N}},\ {\left\Vert{h}\right\Vert}_{W^{m+k,2}(B^{1})}\leq
 c_{k+1}R^{-(k+1)m}{\left\Vert{h}\right\Vert}_{L^{2}(B)}.$\ \par 
\ \par 
The Sobolev embedding theorems, valid in a these balls, give
 that, for any $l\in {\mathbb{N}},\ h\in {\mathcal{C}}^{l}(B^{1}).$
 Then $\displaystyle h\in {\mathcal{C}}^{\infty }(B^{1}).$\ \par 
Because the ${\mathcal{C}}^{\infty }$ regularity is a local property,
 we get that $\displaystyle h\in {\mathcal{C}}^{\infty }(M).$
 $\hfill\blacksquare $\ \par 

\begin{lem}
There is a linear projection from $L^{r}_{G}(M)$ to $\displaystyle
 {\mathcal{H}}_{G}^{2}.$
\end{lem}
\quad Proof.\ \par 
We set\ \par 
\quad \quad \quad \quad \quad $\displaystyle \forall v\in L^{r}_{G}(M),\ H(v):=\sum_{j=1}^{N}{{\left\langle{v,e_{j}}\right\rangle}e_{j}}$\
 \par 
where $\lbrace e_{j}\rbrace _{j=1,...,N}$ is an orthonormal basis
 for $\displaystyle {\mathcal{H}}_{G}^{2}.$ This is meaningful
 because $v\in L^{r}_{G}(M)$ can be integrated against $e_{j}\in
 {\mathcal{H}}_{G}^{2}\subset {\mathcal{C}}^{\infty }(M).$ Moreover
 we have $v-H(v)\in L^{r}_{G}(M)\cap {\mathcal{H}}_{G}^{\perp
 }$ in the sense that $\forall h\in {\mathcal{H}}_{G}^{2},\ {\left\langle{v-H(v),\
 h}\right\rangle}=0;$ it suffices to test on $h:=e_{k}.$ We get\ \par 
\quad \quad \quad $\displaystyle {\left\langle{v-H(v),\ e_{k}}\right\rangle}={\left\langle{v,e_{k}}\right\rangle}-{\left\langle{\sum_{j=1}^{N}{{\left\langle{v,e_{j}}\right\rangle}e_{j},e_{k}}}\right\rangle}={\left\langle{v,e_{k}}\right\rangle}-{\left\langle{v,e_{k}}\right\rangle}=0.$\
 \par 
This ends the proof. $\hfill\blacksquare $\ \par 

\begin{prop}
We have a direct decomposition:\par 
\quad \quad \quad $\displaystyle L^{r}_{G}(M)={\mathcal{H}}_{G}^{2}\oplus \mathrm{I}\mathrm{m}D(W^{2,r}_{G}(M)).$
\end{prop}

      Proof.\ \par 
Let $v\in L^{r}_{G}(M).$ Set $h:=H(v)\in {\mathcal{H}}^{2}_{G},$
 and $\omega :=v-h.$ We have that $\forall k\in {\mathcal{H}}^{2}_{G},\
 {\left\langle{\omega ,k}\right\rangle}={\left\langle{v-H(v),\
 k}\right\rangle}=0.$ Hence we can solve $Du=\omega $ with $u\in
 W^{2,r}_{G}(M)\cap L^{2}_{G}(M).$ So we get $v=h+Du$ which means:\ \par 
\ \par 
\quad \quad \quad $\displaystyle L^{r}_{G}(M)={\mathcal{H}}_{G}^{2}+\mathrm{I}\mathrm{m}D(W^{2,r}_{G}(M)).$\
 \par 
\ \par 
The decomposition is direct because if $\omega \in {\mathcal{H}}^{2}_{G}\cap
 \mathrm{I}\mathrm{m}D(W^{2,r}_{G}(M)),$ then $\omega \in {\mathcal{C}}^{\infty
 }(M)$ and\ \par 
\quad \quad \quad $\omega =Du\Rightarrow \forall k\in {\mathcal{H}}^{2}_{G},\ \omega
 \perp k,$\ \par 
so choosing $k=\omega \in {\mathcal{H}}^{2}_{G}$ we get $\displaystyle
 {\left\langle{\omega ,\omega }\right\rangle}=0$ hence $\omega
 =0.$ The proof is complete. $\hfill\blacksquare $\ \par 
\ \par 
\quad In the special case where $D$ is the Hodge Laplacian, we already
 seen~\cite{RSMcompact17} that we recover this way the strong
 $L^{r}$ Hodge decomposition without using Gaffney's inequalities.\ \par 

\section{Case of compact manifold with a smooth boundary.~\label{LIR34}}
\quad Let $N$ be a ${\mathcal{C}}^{\infty }$ smooth connected riemannian
 manifold  compact with a ${\mathcal{C}}^{\infty }$ smooth boundary
 $\partial N.$ We want to show how the results in case of a compact
 boundary-less manifold apply to this case.\ \par 
\quad First we know that a neighborhood $V$ of $\partial N$ in $N$
 can be seen as $\partial N{\times}\lbrack 0,\delta \rbrack $
 by~\cite[Theorem 5.9 p. 56]{Munkres61}  or by ~\cite[Th\'eor\`eme
 (28) p. 1-21]{Courrege65}. This allows us to "extend" slightly $N:$\ \par 
\quad we have $N=(N\backslash V)\cup V\simeq (N\backslash V)\cup (\partial
 N{\times}\lbrack 0,\delta \rbrack ).$ So we set $M:=(N\backslash
 V)\cup (\partial N{\times}\lbrack 0,\delta +\epsilon \rbrack ).$\ \par 
Then $M$ can be seen as a riemannian manifold with boundary $\partial
 M\simeq \partial N$ and such that $\bar N\subset M.$\ \par 
\quad Now a classical way to get rid of a "annoying boundary" of a
 manifold is to use its "double". For instance: Duff~\cite{Duff52},
 H\"ormander~\cite[p. 257]{Hormand94}. Here we copy the following
 construction from Guneysu and Pigola~\cite[Appendix B]{GuneysuPigola}.\ \par 
\quad The "Riemannian double" $\Gamma :=\Gamma (M)$ of $M,$ obtained
 by gluing two copies, $M$ and $M_{2},$ of $M$ along $\partial
 M,$ is a compact Riemannian manifold without boundary. Moreover,
 by its very construction, it is always possible to assume that
 $\displaystyle \Gamma $ contains an isometric copy of the original
 manifold $N.$ We shall also write $N$ for its isometric copy
 to ease notation.\ \par 
\ \par 
\quad We extend the operator $D$ to $M$ smoothly by extending smoothly
 its coefficients, and because $D$ is strictly elliptic, choosing
 $\epsilon $ small enough, we get that the extension is still
 an elliptic operator on $M.$ Then we take a ${\mathcal{C}}^{\infty
 }$ function $\chi $ with compact support on $M\subset \Gamma
 $ such that: $0\leq \chi \leq 1;\ \chi \equiv 1$ on $N;$ and
 we consider $\tilde D:=\chi D+(1-\chi )D_{2}$ where $D_{2}$
 is the operator $D$ on the copy $M_{2}$ of $M.$ Then $\tilde
 D\equiv D$ on $N$ and is elliptic on $\Gamma .$\ \par 
\quad Now we shall use Definition~\ref{LIR38} from the introduction,
 we recall it here for the reader convenience\ \par 

\begin{defin}
We shall say that $D$ has the weak maximum property, WMP, if,
 for any smooth $DG$-harmonic $h,$ i.e. $G$-form such that $Dh=0$
 in $M,$ smooth up to the boundary $\partial M,$ which is flat
 on $\partial M,$ i.e. zero on $\partial M$ with all its derivatives,
 then $h$ is zero in $M.$
\end{defin}
\quad Of course if there is a maximum principle for $D$ then WMP is
 true. This is the case for smoothly bounded open sets in ${\mathbb{R}}^{n}$
 by a Theorem of S. Agmon~\cite{Agmon60} for functions and by
 ~\cite[Theorem 4.2, p. 59]{AgDougNir64} in the case $G=\Lambda
 ^{p}(M)$ of $p$-forms on $M.$\ \par 
\quad Because this maximum principle is \emph{not} local, I don't know
 what happen on a compact riemannian manifold with smooth boundary
 for general elliptic operator, even in the case $\displaystyle
 G=\Lambda ^{p}(M).$\ \par 
Nevertheless the Hodge laplacian in a riemannian manifold has
 the UCP for $p$-forms by a difficult result  by N. Aronszajn,
 A. Krzywicki and J. Szarski~\cite{Aronszajn62} hence it has the WMP too.\ \par 
\ \par 
\quad The main lemma of this section is:\ \par 

\begin{lem}
~\label{tD0}Let $\omega \in L^{r}_{G}(N),$ then we can extend
 it to $\omega '\in L^{r}_{G}(\Gamma )$ such that: $\forall h\in
 {\mathcal{H}}_{G}(\Gamma ),\ {\left\langle{\omega ',h}\right\rangle}_{\Gamma
 }=0$ provided that the operator $D$ has the WMP for the $D$-harmonic $G$-forms.
\end{lem}
\quad Proof.\ \par 
Recall that ${\mathcal{H}}_{G}(\Gamma ):=\mathrm{k}\mathrm{e}\mathrm{r}D^{*}\cap
 L^{2}_{G}(\Gamma )$ is of finite dimension $K_{G}$  and ${\mathcal{H}}_{G}(\Gamma
 )\subset {\mathcal{C}}^{\infty }(\Gamma )$ by Lemma~\ref{LIR7}.\ \par 
\quad Make an orthonormal basis $\lbrace e_{1},...,e_{K_{G}}\rbrace
 $ of ${\mathcal{H}}_{G}(\Gamma )$ with respect to $L^{2}_{G}(\Gamma
 ),$ by the Gram-Schmidt procedure so ${\left\langle{e_{j},e_{k}}\right\rangle}_{\Gamma
 }:=\int_{\Gamma }{e_{j}e_{k}dv}=\delta _{jk}.$\ \par 
Set $\lambda _{j}:={\left\langle{\omega {\11}_{N},e_{j}}\right\rangle}={\left\langle{\omega
 ,e_{j}{\11}_{N}}\right\rangle},\ j=1,...,K_{G},$ this makes
 sense since $e_{j}\in {\mathcal{C}}^{\infty }(\Gamma )\Rightarrow
 e_{j}\in L^{\infty }(\Gamma ),$ because $\Gamma $ is compact.\ \par 
We shall see that the system $\displaystyle \lbrace e_{k}{\11}_{\Gamma
 \backslash N}\rbrace _{k=1,...,K_{G}}$ is a free one. Suppose
 this is not the case, then it will exist $\gamma _{1},...,\gamma
 _{K_{G}},$ not all zero, such that $\sum_{k=1}^{K_{G}}{\gamma
 _{k}e_{k}{\11}_{\Gamma \backslash N}}=0$ in $\Gamma \backslash
 N.$ But the function $h:=\sum_{k=1}^{K_{G}}{\gamma _{k}e_{k}}$
 is in ${\mathcal{H}}_{G}(\Gamma )$ and $h$ is zero in $\Gamma
 \backslash N$ which is non void, hence $h$ is flat on $\partial
 N.$ Then $h\equiv 0$ in $\Gamma $ by the WMP. But this is not
 possible because the $e_{k}$ make a basis for ${\mathcal{H}}_{G}(\Gamma
 ).$ So the system $\lbrace e_{k}{\11}_{\Gamma \backslash N}\rbrace
 _{k=1,...,K_{G}}$ is a free one.\ \par 
\quad We set $\gamma _{jk}:={\left\langle{e_{k}{\11}_{\Gamma \backslash
 N},e_{j}{\11}_{\Gamma \backslash N}}\right\rangle}$ hence we
 have that $\mathrm{d}\mathrm{e}\mathrm{t}\lbrace \gamma _{jk}\rbrace
 \neq 0.$ So we can solve the linear system to get $\lbrace \mu
 _{k}\rbrace $ such that\ \par 
\quad \quad \quad \begin{equation}  \forall j=1,...,K_{G},\ \sum_{k=1}^{K_{G}}{\mu
 _{k}{\left\langle{e_{k}{\11}_{\Gamma \backslash N},e_{j}}\right\rangle}}=\lambda
 _{j}.\label{RSMH1}\end{equation}\ \par 
We put $\omega '':=\sum_{j=1}^{K_{G}}{\mu _{j}e_{j}{\11}_{\Gamma
 \backslash N}}$ and $\omega ':=\omega {\11}_{N}-\omega ''{\11}_{\Gamma
 \backslash N}=\omega -\omega ''.$ From~(\ref{RSMH1}) we get\ \par 
\quad \quad \quad $\displaystyle \forall j=1,...,K_{G},\ {\left\langle{\omega ',e_{j}}\right\rangle}_{\Gamma
 }={\left\langle{\omega ,e_{j}}\right\rangle}-{\left\langle{\omega
 '',e_{j}}\right\rangle}=\lambda _{j}-\sum_{k=1}^{K_{G}}{\mu
 _{k}{\left\langle{e_{k}{\11}_{\Gamma \backslash N},e_{j}}\right\rangle}}=0.$\
 \par 
So the $G$-form $\omega '$ is orthogonal to ${\mathcal{H}}_{G}.$
 Moreover $\omega '_{\mid N}=\omega $ and clearly $\omega ''\in
 L^{r}_{G}(\Gamma )$ being a finite combination of $\displaystyle
 e_{j}{\11}_{\Gamma \backslash N},$ so $\omega '\in L^{r}_{G}(\Gamma
 )$ because $\omega $ itself is in $\displaystyle L^{r}_{G}(\Gamma
 ).$ The proof is complete. $\hfill\blacksquare $\ \par 
\ \par 
\quad Now let $\omega \in L^{r}_{G}(N)$ and see $N$ as a subset of
 $\Gamma ;$ then extend $\omega $ as $\omega '$ to $\Gamma $
 by Lemma~\ref{tD0}.\ \par 
By the results on the compact manifold $\Gamma ,$ because $\omega
 '\perp {\mathcal{H}}_{G}(\Gamma ),$ we get that there exists
 $u'\in W^{m,r}_{G}(\Gamma ),$ such that $Du'=\omega ';$ hence
 if $u$ is the restriction of $u'$ to $N$ we get  $u\in W^{m,r}_{G}(N),\
 Du=\omega $ in $N.$\ \par 
Hence we proved\ \par 

\begin{thm}
~\label{tD1}Let $N$ be a smooth compact riemannian manifold with
 smooth boundary $\partial N.$ Let $\displaystyle \omega \in
 L^{r}_{G}(N).$ There is a $G$-form $u\in W^{m,r}_{G}(N),$ such
 that $Du=\omega $ and ${\left\Vert{u}\right\Vert}_{W^{m,r}_{G}(N)}\leq
 c{\left\Vert{\omega }\right\Vert}_{L^{r}_{G}(N)},$ provided
 that the operator $D$ has the WMP for the $D$-harmonic $G$-forms.
\end{thm}

\begin{rem}
I had the hope that the WMP condition be also necessary, but
 this is not the case as the Theorem~\ref{mLIR27} shows.
\end{rem}

\section{Relations with the local existence of solutions.~\label{LIR35}}
\quad Let $(M,g)$ be a ${\mathcal{C}}^{\infty }$ smooth compact riemannian
 manifold without boundary.\ \par 
\quad Let $D:G\rightarrow G$ be a linear differential operator of order
 $m$ with ${\mathcal{C}}^{\infty }$ coefficients. \ \par 
As above we suppose that $D$ is elliptic in the sense of Definition~\ref{LIR39}.\
 \par 
\quad Let $x\in M$ and take a ball $B:=B(x,R).$ We suppose that $\omega
 \in L_{G}^{2}(B)$ and we want to solve $Du=\omega .$ For this
 we shall extend $\omega $ as $\omega '\in L^{2}_{G}(M)$ in the
 whole of $M$ with $\omega '\perp {\mathcal{H}}_{G}(M):=\mathrm{k}\mathrm{e}\mathrm{r}D^{*}$
 in order to apply Theorem~\ref{LIR4}.\ \par 
\quad Consider $\omega :=\omega {\11}_{B}$ the trivial extension of
 $\omega $ to $M.$ We have, with $P_{h}$ the orthogonal projection
 on ${\mathcal{H}}_{G}(M),\ h:=P_{h}\omega .$ Set $N:=K_{G}$
 the finite dimension of $\displaystyle {\mathcal{H}}_{G}(M).$
 Take an orthonormal basis $\lbrace e_{1},...,\ e_{N}\rbrace
 $ of ${\mathcal{H}}_{G}(M),$ then we have\ \par 
\quad \quad \quad $h:=\sum_{j=1}^{N}{h_{j}e_{j}}.$\ \par 
If $h=0,$ we set $\omega '=\omega $ and we are done. If not let
 the radius $R$ of the ball $B$ be small enough to have\ \par 
\ \par 
\quad $\displaystyle {\left\Vert{e_{1}{\11}_{B}}\right\Vert}\leq \frac{1}{4{\sqrt{N}}},\
 ...,\ {\left\Vert{e_{N}{\11}_{B}}\right\Vert}\leq \frac{1}{4{\sqrt{N}}}.$\
 \par 
\ \par 
This is possible because the $e_{j}$ are in ${\mathcal{C}}^{\infty
 }(M)$ so if $B$ is small enough we have $\displaystyle {\left\Vert{e_{j}{\11}_{B}}\right\Vert}\leq
 \frac{1}{4{\sqrt{N}}},$ and we have a finite number of such conditions.\ \par 
\quad We set $\displaystyle \omega _{1}:={\11}_{B^{c}}\sum_{j=1}^{N}{h_{j}e_{j}}.$
 Then\ \par 
\quad \quad \quad $\displaystyle {\left\Vert{\omega _{1}}\right\Vert}^{2}:=\int_{B^{c}}{\left\vert{\sum_{j=1}^{N}{h_{j}e_{j}}}\right\vert
 ^{2}dv\leq }\int_{M}{\left\vert{\sum_{j=1}^{N}{h_{j}e_{j}}}\right\vert
 ^{2}dv\leq }{\left\Vert{h}\right\Vert}^{2}.$\ \par 
And\ \par 
\quad \quad \quad $\displaystyle {\left\Vert{h-\omega _{1}}\right\Vert}\leq \sum_{j=1}^{N}{\left\vert{h_{j}}\right\vert
 {\left\Vert{{\11}_{B}e_{j}}\right\Vert}}\leq \frac{1}{4}\sum_{j=1}^{N}{\left\vert{h_{j}}\right\vert
 }\leq {\sqrt{N}}{\left\Vert{h}\right\Vert}\frac{1}{4{\sqrt{N}}}=\frac{1}{4}{\left\Vert{h}\right\Vert}.$\
 \par 
\ \par 
Hence, because $P_{h}$ has norm one,\ \par 
\ \par 
\quad \quad \quad $\displaystyle {\left\Vert{h-P_{h}\omega _{1}}\right\Vert}={\left\Vert{P_{h}h-P_{h}\omega
 _{1}}\right\Vert}\leq {\left\Vert{h-\omega _{1}}\right\Vert}\leq
 \frac{1}{4}{\left\Vert{h}\right\Vert}.$\ \par 
\ \par 
Now we set $h_{1}:=h-P_{h}\omega _{1}.$ Then ${\left\Vert{h_{1}}\right\Vert}\leq
 \frac{1}{4}{\left\Vert{h}\right\Vert}$ and we have $\displaystyle
 h_{1}:=\sum_{j=1}^{N}{h^{1}_{j}e_{j}}.$ So we set $\omega _{2}:={\11}_{B^{c}}\sum_{j=1}^{N}{h^{1}_{j}e_{j}}.$
 We have the same way:\ \par 
\quad \quad \quad $\displaystyle {\left\Vert{\omega _{2}}\right\Vert}\leq {\left\Vert{h_{1}}\right\Vert}\leq
 \frac{1}{4}{\left\Vert{h}\right\Vert}$ and $\displaystyle {\left\Vert{h_{1}-P_{h}\omega
 _{2}}\right\Vert}\leq \frac{1}{4}{\left\Vert{h_{1}}\right\Vert}\leq
 \frac{1}{4^{2}}{\left\Vert{h}\right\Vert}.$\ \par 
At the step $k$ we get:\ \par 
\quad \quad \quad $\displaystyle {\left\Vert{h_{k}-P_{h}\omega _{k+1}}\right\Vert}\leq
 \frac{1}{4}{\left\Vert{h_{k}}\right\Vert}\leq \frac{1}{4^{k}}{\left\Vert{h}\right\Vert}$
 and $\displaystyle {\left\Vert{\omega _{k+1}}\right\Vert}\leq
 \frac{1}{4^{k}}{\left\Vert{h}\right\Vert}.$\ \par 
We set $\omega '':=\sum_{j=1}^{\infty }{\omega _{j}}.$ We get
 that the series converges in norm $L^{2}(M)$ and $P_{h}\omega ''=h.$\ \par 
Setting $\omega ':=\omega -\omega '',$ we get that $\omega '=\omega
 $ on $B$ and $P_{h}(\omega ')=0,$ which means that $\omega '\perp
 {\mathcal{H}}_{G}(M).$\ \par 
\quad We can apply Theorem~\ref{LIR4} to get $Du'=\omega '$ with $\displaystyle
 u'\in L^{2}_{G}(M)$ because $\omega '\perp {\mathcal{H}}_{G}.$
 We set $u:=u'_{\mid B}$ in $B$ to have $Du=\omega $ in $B.$\ \par 
So we proved:\ \par 

\begin{thm}
~\label{r1} Let $x$ in $M.$ There is a $R_{0}(x)>0$ such that
 for any $0<R\leq R_{0}$ if $\displaystyle \omega \in L_{G}^{2}(B)$
 with $B:=B(x,R)$ there is a $\displaystyle u\in L_{G}^{2}(B)$
 such that $Du=\omega $ and ${\left\Vert{u}\right\Vert}_{L_{G}^{2}(B)}\lesssim
 {\left\Vert{\omega }\right\Vert}_{L_{G}^{2}(B)}.$
\end{thm}
\quad To get the $\displaystyle L_{G}^{r}(B)$ case for $r>2,$ we proceed
 as in the proof of Theorem~\ref{LIR1}.\ \par 

\begin{thm}
~\label{mLIR27}Under the assumptions above, for any $x\in M$
 and $r\geq 2,$ there is a positive constant $\displaystyle c_{f}$
 such that, if $\omega \in L^{r}(B),$ there is a$\ u\in L^{t}(B^{1})$
 with $\displaystyle \frac{1}{t}=\frac{1}{r}-\tau ,$ such that
 $Du=\omega $ and ${\left\Vert{u}\right\Vert}_{L^{t}(B^{1})}\leq
 c_{f}{\left\Vert{\omega }\right\Vert}_{L^{r}(B)}.$\par 
\quad Moreover we have $u\in W_{G}^{m,r}(B^{1})$ with control of the norm.
\end{thm}
\quad Proof.\ \par 
Let $r\geq 2$ and $\displaystyle \omega \in L^{r}_{G}(B).$ Because
 $B$ is relatively compact and $dv$ is $\sigma $-finite, we have
 that $\displaystyle \omega \in L^{2}_{G}(B).$ Theorem~\ref{r1}
 gives that there is a $u\in L^{2}_{G}(B)$ such that $Du=\omega
 .$ Now we proceed exactly as in the proof of Theorem~\ref{LIR1},
 using the same induction procedure. $\hfill\blacksquare $\ \par 
\ \par 
\quad So we proved the local existence of solutions with estimates;
 this is an already known theorem in ${\mathbb{R}}^{n},$ hence
 also locally in $M$ (see for instance~\cite{DongKim11}).   
   This means also that the LIR condition is stronger than the
 local existence of solutions with estimates. These solutions
 were the basis of the Raising Steps Method, see~\cite{HodgeNonCompact15}.\
 \par 

\section{The non compact case.~\label{LIR36}}
\quad We shall use the same ideas as in~\cite{HodgeNonCompact15} to
 go from the compact case to the non compact one.\ \par 
\quad In order to deal with $G$-forms in the non compact case, we have
 to warranty that the bundle $G$ has trivializing charts defined
 on balls of the covering ${\mathcal{C}}_{\epsilon }.$\ \par 

\begin{defin}
We say that the bundle $\displaystyle G:=(H,\pi ,M)$ is compatible
 with the covering ${\mathcal{C}}_{\epsilon }$ if there is a
 $\epsilon >0$ such that, for any ball $B\in {\mathcal{C}}_{\epsilon
 },$ the chart $(B,\varphi )$ is a trivializing map of the bundle
 $G.$ Precisely this means that $G\simeq \varphi (B){\times}{\mathbb{R}}^{N}$
 where $N$ is the dimension of $H$ and the equivalence has bounds
 independent of $\displaystyle B\in {\mathcal{C}}_{\epsilon }.$
\end{defin}

\begin{exa}
The bundle of $p$-forms in a riemannian manifold $(M,g)$ is compatible.
 To see this take a ball $B(x,R)\in {\mathcal{C}}_{\epsilon },$
 then we have that $\displaystyle (1-\epsilon )\delta _{ij}\leq
 g_{ij}\leq (1+\epsilon )\delta _{ij}$ in $\displaystyle B(x,R)$
 as bilinear forms, so, because $\epsilon <1,$ the $1$-forms
 $dx_{j},\ j=1,...,n$ are "almost" orthonormal and hence linearly
 independent. This gives that the co-tangent bundle $T^{*}M$
 is equivalent to $T^{*}{\mathbb{R}}^{n}$ over $B,$ the constants
 depending only on $\epsilon .$\par 
\quad By tensorisation we get the same for the bundle of  $p$-forms.
\end{exa}
\quad From now on we shall always suppose that the bundle $\displaystyle
 G:=(H,\pi ,M)$ is compatible with the covering ${\mathcal{C}}_{\epsilon
 }.$\ \par 
\ \par 
In subsection~\ref{CL24} we define a Vitali type covering ${\mathcal{C}}_{\epsilon
 }$ by balls suited to our "admissible balls" (see Definition~\ref{mLIR25}).
 We use these notions now.\ \par 

\begin{defin}
~\label{mLIR26}We shall say that the hypothesis (UEAB) is fulfilled
 for the operator $D$ if $D$ has smooth ${\mathcal{C}}^{1}(M)$
 coefficients.\par 
Moreover we ask that $D$ be uniformly elliptic as in Definition~\ref{LIR39}.
\end{defin}
\quad We start with $\omega \in L^{2}_{G}(M),$ by the (THL2p) hypothesis,
 provided that $\omega \perp \mathrm{k}\mathrm{e}\mathrm{r}D^{*},$
 there is a $G$-form $u\in L^{2}_{G}(M)$ such that $Du=\omega
 .$ Moreover, because $\displaystyle L^{2}_{G}(M)$ is a Hilbert
 space, the $u\in L^{2}_{G}(M),\ Du=\omega $ with the smallest
 norm, is given linearly with respect to $\omega .$ This means
 that we have a bounded linear operator $\displaystyle S:\ L^{2}_{G}(M)\rightarrow
 L^{2}_{G}(M)$ such that $D(S\omega )=\omega $ provided that
 $\displaystyle \omega \perp \mathrm{k}\mathrm{e}\mathrm{r}D^{*}.$\ \par 
\ \par 
\quad The local\emph{ elliptic inequalities} by Theorem~\ref{LIR0}
 become uniform by the hypothesis (UEAB):\ \par 

\begin{cor}
Let $D$ be an operator of order $m$ acting on sections of $G$
 in the complete riemannian manifold $M.$ Suppose that $D$ verifies
 (UEAB). Then, for any $B_{x}:=B(x,R)\in {\mathcal{C}}_{\epsilon
 }$ and $\displaystyle B^{1}_{x}:=B(x,R/2),$ we have, with $D$
 with ${\mathcal{C}}^{1}(M)$ coefficients:\par 
\quad \quad \quad $\displaystyle {\left\Vert{u}\right\Vert}_{W^{m,r}_{G}(B^{1}_{x})}\leq
 c_{1}{\left\Vert{Du}\right\Vert}_{L^{r}_{G}(B_{x})}+c_{2}R^{-m}{\left\Vert{u}\right\Vert}_{L^{r}_{G}(B_{x})}.$\par
 
The hypotheses (UEAB) are precisely done to warranty  that the
 constants $c_{1},c_{2}$ depend only on $n=\mathrm{d}\mathrm{i}\mathrm{m}_{{\mathbb{R}}}M,\
 r$ and $\epsilon .$
\end{cor}
\quad With $\displaystyle t=S_{m}(r),$ we get, by Lemma~\ref{3S2} from
 the Appendix,\ \par 
\quad \quad \quad $\displaystyle {\left\Vert{u}\right\Vert}_{L^{t}_{G}(B(x,R))}\leq
 CR^{-m}\ {\left\Vert{u}\right\Vert}_{W^{m,r}_{G}(B(x,R))}.$\ \par 
When there is no ambiguity we shall omit the subscript $G,$ i.e.
 $\displaystyle L^{2}_{G}(B)$ becomes $\displaystyle L^{2}(B),$ etc...\ \par 

\begin{lem}
~\label{mLIR10}We have, with $B^{l}:=B(x,2^{-l}R)$ and $t_{0}=2,\
 B^{0}=B(x,R),$ the a priori estimates:\par 
\quad \quad \quad $\displaystyle R^{(l+1)m}{\left\Vert{u}\right\Vert}_{L^{t_{l}}(B^{l})}\leq
 \sum_{j=1}^{l}{c_{j}R^{(l-j+1)m}{\left\Vert{Du}\right\Vert}_{L^{l-j}(B^{l-j})}}+c_{l+1}{\left\Vert{u}\right\Vert}_{L^{2}(B)}.$\par
 
And\par 
\quad \quad \quad $\displaystyle R^{(l+2)m}{\left\Vert{u}\right\Vert}_{W^{m,t_{l}}(B^{l+1})}\leq
 c_{0}R^{(l+2)m}{\left\Vert{Du}\right\Vert}_{L^{t_{l}}(B^{l})}+\sum_{j=1}^{l}{c_{j}R^{(l-j+1)m}{\left\Vert{Du}\right\Vert}_{L^{t_{l-j}}(B^{l-j})}}+c_{l+1}{\left\Vert{u}\right\Vert}_{L^{2}(B)}.$
\end{lem}

      Proof.\ \par 
From the LIR, Theorem~\ref{LIR0}, we have\ \par 
\ \par 
\quad \quad \quad $\displaystyle \forall B\in {\mathcal{C}}_{\epsilon },\ {\left\Vert{u}\right\Vert}_{W^{m,2}(B^{1})}\leq
 c_{1}{\left\Vert{D(u)}\right\Vert}_{L^{2}(B)}+c_{2}R^{-m}{\left\Vert{u}\right\Vert}_{L^{2}(B)}.$\
 \par 
\ \par 
\quad Now we shall use the local Sobolev embedding theorem, Lemma~\ref{3S2},
 to get:\ \par 
\ \par 
\quad \quad \quad $\displaystyle \forall B\in {\mathcal{C}}_{\epsilon },\ {\left\Vert{u}\right\Vert}_{L^{t_{1}}(B^{1})}\leq
 CR^{-m}{\left\Vert{u}\right\Vert}_{W^{m,2}(B)}$\ \par 
so we get\ \par 
\quad \quad $\forall B\in {\mathcal{C}}_{\epsilon },\ {\left\Vert{u}\right\Vert}_{L^{t_{1}}(B^{1})}\leq
 c_{1}R^{-m}{\left\Vert{Du}\right\Vert}_{L^{2}(B)}+c_{2}R^{-2m}{\left\Vert{u}\right\Vert}_{L^{2}(B)}$\
 \par 
\ \par 
with $\frac{1}{t_{1}}:=\frac{1}{2}-\frac{m}{n}\iff t_{1}:=S_{m}(2).$\ \par 
\ \par 
\quad $\bullet $ If $t_{1}\geq r,$ then we get still by the LIR, Theorem~\ref{LIR0}:\
 \par 
\quad \quad \quad \begin{equation}  \forall B\in {\mathcal{C}}_{\epsilon },\ {\left\Vert{u}\right\Vert}_{W^{m,t_{1}}(B^{2})}\leq
 c_{1}{\left\Vert{Du}\right\Vert}_{L^{t_{1}}(B^{1})}+c_{2}R^{-m}{\left\Vert{u}\right\Vert}_{L^{t_{1}}(B^{1})}.\label{nC9}\end{equation}\
 \par 
Putting the estimate of $\displaystyle {\left\Vert{u}\right\Vert}_{L^{t_{1}}(B^{1})}$
 in~(\ref{nC9}) we get\ \par 
\ \par 
\quad \quad \quad $\displaystyle {\left\Vert{u}\right\Vert}_{W^{m,t_{1}}(B^{2})}\leq
 c_{1}{\left\Vert{Du}\right\Vert}_{L^{t_{1}}(B^{1})}+c_{2}R^{-m}{\left({c_{1}{\left\Vert{Du}\right\Vert}_{L^{2}(B)}+c_{2}R^{-m}{\left\Vert{u}\right\Vert}_{L^{2}(B)}}\right)}$\
 \par 
\ \par 
so, with suitable constants\ \par 
\ \par 
\quad \quad \quad $\displaystyle {\left\Vert{u}\right\Vert}_{W^{m,t_{1}}(B^{2})}\leq
 c_{1}{\left\Vert{Du}\right\Vert}_{L^{t_{1}}(B^{1})}+c_{2}R^{-m}{\left\Vert{Du}\right\Vert}_{L^{2}(B)}+c_{3}R^{-2m}{\left\Vert{u}\right\Vert}_{L^{2}(B)}.$\
 \par 
\ \par 
Putting the powers of $R$ on the other side to isolate $\displaystyle
 {\left\Vert{u}\right\Vert}_{L^{2}(B)},$ we get\ \par 
\ \par 
\quad \quad \quad $\displaystyle R^{2m}{\left\Vert{u}\right\Vert}_{W^{m,t_{1}}(B^{2})}\leq
 c_{1}R^{2m}{\left\Vert{Du}\right\Vert}_{L^{t_{1}}(B^{1})}+c_{2}R^{m}{\left\Vert{Du}\right\Vert}_{L^{2}(B)}+c_{3}{\left\Vert{u}\right\Vert}_{L^{2}(B)}.$\
 \par 
\ \par 
We iterate, using again the local Sobolev embedding theorem,
 Lemma~\ref{3S2},\ \par 
\ \par 
\quad \quad $\displaystyle u\in L^{t_{2}}(B^{2}),\ {\left\Vert{u}\right\Vert}_{L^{t_{2}}(B^{2})}\leq
 cR^{-m}{\left\Vert{u}\right\Vert}_{W^{m,t_{1}}(B^{2})}$\ \par 
hence\ \par 
\quad \quad $\displaystyle R^{3m}{\left\Vert{u}\right\Vert}_{L^{t_{2}}(B^{2})}\leq
 c_{1}R^{2m}{\left\Vert{Du}\right\Vert}_{L^{t_{1}}(B^{1})}+c_{2}R^{m}{\left\Vert{Du}\right\Vert}_{L^{2}(B)}+c_{3}{\left\Vert{u}\right\Vert}_{L^{2}(B)}.$\
 \par 
\ \par 
with $\displaystyle \frac{1}{t_{2}}:=\frac{1}{t_{1}}-\frac{m}{n}=\frac{1}{2}-\frac{2m}{n}\iff
 t_{2}:=S_{2m}(2).$ The LIR gives again:\ \par 
\ \par 
\quad \quad \quad $\displaystyle {\left\Vert{u}\right\Vert}_{W^{m,t_{2}}(B^{3})}\leq
 c_{1}{\left\Vert{Du}\right\Vert}_{L^{t_{2}}(B^{2})}+c_{2}R^{-m}{\left\Vert{u}\right\Vert}_{L^{t_{2}}(B^{2})}$\
 \par 
so\ \par 
\quad \quad \quad $\displaystyle R^{4m}{\left\Vert{u}\right\Vert}_{W^{m,t_{2}}(B^{3})}\leq
 c_{1}R^{4m}{\left\Vert{Du}\right\Vert}_{L^{t_{2}}(B^{2})}+c_{2}R^{3m}{\left\Vert{u}\right\Vert}_{L^{t_{2}}(B^{2})}$\
 \par 
hence\ \par 
\quad \quad \quad $\displaystyle R^{4m}{\left\Vert{u}\right\Vert}_{W^{m,t_{2}}(B^{3})}\leq
 c_{1}R^{4m}{\left\Vert{Du}\right\Vert}_{L^{t_{2}}(B^{2})}+c_{2}R^{2m}{\left\Vert{Du}\right\Vert}_{L^{t_{1}}(B^{1})}+c_{3}R^{m}{\left\Vert{Du}\right\Vert}_{L^{2}(B)}+c_{4}{\left\Vert{u}\right\Vert}_{L^{2}(B)}.$\
 \par 
\ \par 
\quad Iterating the same way we get\ \par 
\ \par 
\quad \quad \quad $\displaystyle R^{(l+1)m}{\left\Vert{u}\right\Vert}_{L^{t_{l}}(B^{l})}\leq
 c_{1}R^{lm}{\left\Vert{Du}\right\Vert}_{L^{t_{l-1}}(B^{l-1})}+c_{2}R^{(l-1)m}{\left\Vert{Du}\right\Vert}_{L^{t_{(l-2)}}(B^{(l-2)})}+\cdot
 \cdot \cdot +$\ \par 
\quad \quad \quad \quad \quad \quad \quad \quad \quad \quad \quad \quad \quad \quad \quad $\displaystyle +c_{l}R^{m}{\left\Vert{Du}\right\Vert}_{L^{2}(B)}+c_{l+1}{\left\Vert{u}\right\Vert}_{L^{2}(B)}.$\
 \par 
\ \par 
Which gives, using the LIR,\ \par 
\ \par 
\quad \quad \quad $\displaystyle {\left\Vert{u}\right\Vert}_{W^{m,t_{l}}(B^{l+1})}\leq
 c_{1}{\left\Vert{Du}\right\Vert}_{L^{t_{l}}(B^{l})}+c_{2}R^{-m}{\left\Vert{u}\right\Vert}_{L^{t_{l}}(B^{l})}$\
 \par 
so\ \par 
\quad \quad \quad $\displaystyle R^{(l+2)m}{\left\Vert{u}\right\Vert}_{W^{m,t_{l}}(B^{l+1})}\leq
 c_{1}R^{(l+2)m}{\left\Vert{Du}\right\Vert}_{L^{t_{l}}(B^{l})}+c_{2}R^{(l+1)m}{\left\Vert{u}\right\Vert}_{L^{t_{l}}(B^{l})}$\
 \par 
and\ \par 
\quad \quad \quad $\displaystyle R^{(l+2)m}{\left\Vert{u}\right\Vert}_{W^{m,t_{l}}(B^{l+1})}\leq
 c_{1}R^{(l+2)m}{\left\Vert{Du}\right\Vert}_{L^{t_{l}}(B^{l})}+c_{2}R^{lm}{\left\Vert{Du}\right\Vert}_{L^{t_{(l-1)}}(B^{(l-1)})}+$\
 \par 
\quad \quad \quad \quad \quad \quad \quad \quad \quad \quad \quad \quad \quad \quad \quad \quad $\displaystyle +c_{3}R^{(l-1)m}{\left\Vert{Du}\right\Vert}_{L^{t_{(l-2)}}(B^{(l-2)})}+\cdot
 \cdot \cdot +c_{l}R^{m}{\left\Vert{Du}\right\Vert}_{L^{2}(B)}+c_{l+1}{\left\Vert{u}\right\Vert}_{L^{2}(B)}.$\
 \par 
Which proves the lemma. $\hfill\blacksquare $\ \par 

\begin{lem}
~\label{mLIR21} We have for $r<t,\ B:=B(x,R),$\par 
\quad \quad \quad $\displaystyle \forall f\in L^{r}(B),\ {\left\Vert{f}\right\Vert}_{L^{r}(B)}\leq
 R^{\frac{1}{r}-\frac{1}{t}}{\left\Vert{f}\right\Vert}_{L^{t}(B)}.$
\end{lem}
\quad Proof.\ \par 
Because the measure $\displaystyle d\mu (x):=\frac{{\11}_{B}(x)}{\left\vert{B}\right\vert
 }dm(x)$ is a probability measure, using that $r<t,$ we have
 ${\left\Vert{f}\right\Vert}_{L^{r}(\mu )}\leq {\left\Vert{f}\right\Vert}_{L^{t}(\mu
 )}$ which implies readily the lemma. $\hfill\blacksquare $\ \par 

\begin{cor}
~\label{mLIR24}Let $\forall j\in {\mathbb{N}},\ \frac{1}{t_{j}}=\frac{1}{2}-\frac{jm}{n}.$
 Fix $r\geq 2,$ we have, for $t_{l-1}<r<t_{l},$\par 
\quad \quad \quad $\displaystyle R^{(\frac{1}{t_{l}}-\frac{1}{r})+(l+1)m}{\left\Vert{u}\right\Vert}_{L^{r}(B^{l})}\leq
 \sum_{j=1}^{l}{c_{j}R^{(l-j+1)m}{\left\Vert{Du}\right\Vert}_{L^{l-j}(B^{l-j})}}+c_{l+1}{\left\Vert{u}\right\Vert}_{L^{2}(B)}.$
\end{cor}

      Proof.\ \par 
By Lemma~\ref{mLIR21} we get $\displaystyle {\left\Vert{u}\right\Vert}_{L^{r}(B^{l})}\leq
 R^{\frac{1}{r}-\frac{1}{t_{l}}}{\left\Vert{u}\right\Vert}_{L^{t}(B^{l})}$
 so by Lemma~\ref{mLIR10} we have\ \par 
\quad \quad \quad $\displaystyle R^{(l+1)m}{\left\Vert{u}\right\Vert}_{L^{r}(B^{l})}\leq
 R^{\frac{1}{r}-\frac{1}{t_{l}}}{\left\Vert{u}\right\Vert}_{L^{t}(B^{l})}\leq
 R^{\frac{1}{r}-\frac{1}{t_{l}}}\sum_{j=1}^{l}{c_{j}R^{(l-j+1)m}{\left\Vert{Du}\right\Vert}_{L^{l-j}(B^{l-j})}}+c_{l+1}R^{\frac{1}{r}-\frac{1}{t_{l}}}{\left\Vert{u}\right\Vert}_{L^{2}(B)}.$\
 \par 
Isolating $\displaystyle {\left\Vert{u}\right\Vert}_{L^{2}(B)}$ we get\ \par 
\quad \quad \quad $\displaystyle R^{(\frac{1}{t_{l}}-\frac{1}{r})+(l+1)m}{\left\Vert{u}\right\Vert}_{L^{r}(B^{l})}\leq
 \sum_{j=1}^{l}{c_{j}R^{(l-j+1)m}{\left\Vert{Du}\right\Vert}_{L^{l-j}(B^{l-j})}}+c_{l+1}{\left\Vert{u}\right\Vert}_{L^{2}(B)}.$\
 \par 
Now we have a finite number of terms, so changing the values
 of the constants, we get\ \par 
\quad \quad \quad $\displaystyle R^{(\frac{r}{t_{l}}-1)+(l+1)mr}{\left\Vert{u}\right\Vert}_{L^{r}(B^{l})}^{r}\leq
 \sum_{j=1}^{l}{c_{j}R^{(l-j+1)mr}{\left\Vert{Du}\right\Vert}_{L^{l-j}(B^{l-j})}^{r}}+c_{l+1}{\left\Vert{u}\right\Vert}_{L^{2}(B)}^{r}.$\
 \par 
which ends the proof of the corollary. $\hfill\blacksquare $\ \par 
\ \par 
\quad We shall use the following weights, with $t_{j}:=S_{jm}(2)\ i.e.\
 \frac{1}{t_{j}}=\frac{1}{2}-\frac{jm}{n}:$\ \par 
\quad \quad \quad $\displaystyle t_{l-1}<r<t_{l},\ v_{r}(x):=R(x)^{(\frac{1}{t_{l}}-\frac{1}{r})+(l+1)m},\
 w_{j}(x)=R^{(l+1-j)m}$\ \par 
and we set:\ \par 
\quad \quad \quad $\displaystyle {\left\Vert{\omega }\right\Vert}_{L^{t_{l-j}}(M,w_{j}^{t_{l-j}})}^{t_{l-j}}:=\int_{M}{\left\vert{\omega
 (x)}\right\vert ^{t_{l-j}}w_{j}(x)^{t_{l-j}}dv(x)}.$\ \par 

\begin{thm}
Under hypotheses (THL2G) and (UEAB), with the weights defined
 above, we have, provided that $\displaystyle \omega \perp \mathrm{k}\mathrm{e}\mathrm{r}D^{*},$
 that there is a $u:=S\omega $ linearly given from $\omega $
 such that $\displaystyle Du=\omega $ and:\par 
\quad \quad \quad $\displaystyle {\left\Vert{u}\right\Vert}_{L^{r}_{G}(M,v_{r}^{r})}\leq
 \sum_{j=1}^{l}{c_{j}{\left\Vert{\omega }\right\Vert}_{L^{t_{l-j}}_{G}(M,w_{j}^{t_{l-j}})}}+c_{l+1}{\left\Vert{\omega
 }\right\Vert}_{L^{2}_{G}(M)}.$
\end{thm}
\quad Proof.\ \par 
By hypothesis (THL2G) for $\omega \in L^{2}_{G}(M)$ with $\displaystyle
 \omega \perp \mathrm{k}\mathrm{e}\mathrm{r}D^{*}$ we set $u:=S\omega
 \in L^{2}_{G}(M).$\ \par 
We have, with hypothesis (UEAB) and using the covering of $M$
 by the $B^{l},$ hence a fortiori by the $B^{j},\ j<l,$\ \par 
\quad \quad \quad \begin{equation} {\left\Vert{ u}\right\Vert}_{L^{r}(M,v_{r}^{r})}^{r}\leq
 \sum_{B\in {\mathcal{C}}_{\epsilon }}{R^{(\frac{1}{t_{l}}-\frac{1}{r})+(l+1)mr}{\left\Vert{u}\right\Vert}_{L^{r}(B^{l})}^{r}}.\label{mLIR23}\end{equation}\
 \par 
Using that the overlap of the covering is bounded by $T,$\ \par 
\quad \quad \quad \begin{equation} \sum_{ B\in {\mathcal{C}}_{\epsilon }}{R^{(l+1-j)mt_{l-j}}{\left\Vert{\omega
 }\right\Vert}_{L^{t_{(l-j)}}(B^{l-j})}^{t_{l-j}}}\leq T{\left\Vert{\omega
 }\right\Vert}_{L^{t_{l-j}}(M,w_{j}^{t_{l-j}})}^{t_{l-j}}.\label{mLIR22}\end{equation}\
 \par 
with $\displaystyle w_{j}(x)=w_{j,l}(x)=R^{(l+1-j)m},$ and for
 any $\gamma ,\ {\left\Vert{\gamma }\right\Vert}_{L^{s}(M,w_{k}^{s})}^{s}:=\int_{M}{\left\vert{\gamma
 (x)w_{k}(x)}\right\vert ^{s}dv(x)}.$\ \par 
Now if $r\geq t_{l-1}\geq t_{l-j},\ j\leq l-1,$ we have $\sum_{j\in
 {\mathbb{N}}}{a_{j}^{r}}\leq {\left({\sum_{j\in {\mathbb{N}}}{a_{j}^{t_{l-j}}}}\right)}^{r/t_{l-j}},$
 so\ \par 
\ \par 
\quad \quad \quad $\displaystyle \sum_{B\in {\mathcal{C}}_{\epsilon }}{R^{(l+1-j)mr}{\left\Vert{\omega
 }\right\Vert}_{L^{t_{(l-j)}}(B^{l-j})}^{r}}\leq {\left({\sum_{B\in
 {\mathcal{C}}_{\epsilon }}{R^{(l+1-j)mt_{l-j}}{\left\Vert{\omega
 }\right\Vert}_{L^{t_{(l-j)}}(B^{l-j})}^{t_{l-j}}}}\right)}^{r/t_{l-j}}.$\ \par 
\ \par 
Using~(\ref{mLIR22}) we get\ \par 
\ \par 
\quad \quad \quad $\displaystyle \sum_{B\in {\mathcal{C}}_{\epsilon }}{R^{(l+1-j)mr}{\left\Vert{\omega
 }\right\Vert}_{L^{t_{(l-j)}}(B^{l-j})}^{r}}\leq T^{r/t_{l-j}}{\left\Vert{\omega
 }\right\Vert}_{L^{t_{l-j}}(M,w_{j}^{t_{l-j}})}^{r}.$\ \par 
\ \par 
Grouping with~(\ref{mLIR23}) we deduce\ \par 
\ \par 
\quad \quad \quad $\displaystyle {\left\Vert{u}\right\Vert}_{L^{r}(M,v_{r}^{r})}^{r}\leq
 \sum_{j=1}^{l}{c_{j}T^{r/t_{l-j}}{\left\Vert{\omega }\right\Vert}_{L^{t_{l-j}}(M,w_{j}^{t_{l-j}})}^{r}}+c_{l+1}{\left\Vert{u}\right\Vert}_{L^{2}(M)}^{r}.$\
 \par 
\ \par 
Changing the constants, we take the $r$ root to get, using the
 hypothesis (THL2G), which says also that $\displaystyle {\left\Vert{u}\right\Vert}_{L^{2}(M)}\leq
 c{\left\Vert{\omega }\right\Vert}_{L^{2}(M)},$\ \par 
\ \par 
\quad \quad \quad $\displaystyle {\left\Vert{u}\right\Vert}_{L^{r}(M,v_{r}^{r})}\leq
 \sum_{j=1}^{l}{c_{j}{\left\Vert{\omega }\right\Vert}_{L^{t_{l-j}}(M,w_{j}^{t_{l-j}})}}+c_{l+1}{\left\Vert{\omega
 }\right\Vert}_{L^{2}(M)}.$\ \par 
\ \par 
The proof is complete. $\hfill\blacksquare $\ \par 

\begin{lem}
~\label{mLIR12}Provided that $\displaystyle \omega \in L^{2}(M)\cap
 L^{t_{k}}(M,R(x)^{\alpha _{k}}),$ with\par 
\quad \quad \quad $\displaystyle \alpha _{j}:=\frac{k+1}{k}m{\times}jt_{j},\ \beta
 _{j}:=(j+1)m{\times}t_{j},$\par 
we have:\par 
\quad \quad \quad $\displaystyle \forall j\leq k,\ \omega \in L^{t_{j}}(M,R^{\beta
 _{j}}),\ {\left\Vert{\omega }\right\Vert}_{L^{t_{j}}(M,R^{\beta
 _{j}})}\leq C\max ({\left\Vert{\omega }\right\Vert}_{L^{t_{k}}(M,R^{\alpha
 _{k}})},\ {\left\Vert{\omega }\right\Vert}_{L^{2}(M)}).$
\end{lem}
\quad Proof.\ \par 
Recall the Stein-Weiss interpolation Theorem ~\cite[Theorem 5.5.1,
 p. 110]{BerghLofstrom76}\ \par 
\quad \quad \quad ${\left({L^{s_{0}}(v_{0}),L^{s_{1}}(v_{1})}\right)}_{\theta ,t}=L^{s}(v),\
 0<\theta <1$ where $v:=v_{0}^{s(1-\theta )/s_{0}}v_{1}^{s\theta
 /s_{1}},\ \frac{1}{s}=\frac{1-\theta }{s_{0}}+\frac{\theta }{s_{1}}.$\ \par 
We choose $s_{0}=2,\ v_{0}=1\ ;\ s_{1}=t_{k}=S_{km}(2),\ s=t_{j}=S_{jm}(2),$
 so $\frac{1}{t_{k}}=\frac{1}{2}-\frac{km}{n},\ \frac{1}{t_{j}}=\frac{1}{2}-\frac{jm}{n}.$
 This fixes $\theta $:\ \par 
\quad \quad \quad $\displaystyle \frac{1}{s}=\frac{1}{t_{j}}=\frac{1}{2}-\frac{jm}{n}=(1-\theta
 )\frac{1}{2}+\theta (\frac{1}{2}-\frac{km}{n})\Rightarrow \theta
 =\frac{j}{k}.$\ \par 
Replacing  $v_{0}=w_{1}^{2}=1,\ v_{1}=w_{2}^{s_{1}}=R(x)^{(k+1)m{\times}t_{k}}$
 and using $\displaystyle v:=v_{0}^{s(1-\theta )/s_{0}}v_{1}^{s\theta
 /s_{1}}$ we get\ \par 
\ \par 
\quad \quad \quad $v=v_{1}^{\frac{s}{s_{1}}{\times}\frac{j}{k}}\Rightarrow \frac{s}{s_{1}}{\times}\frac{j}{k}=\frac{t_{j}}{t_{k}}{\times}\frac{j}{k}\Rightarrow
 v=R(x)^{(k+1)m{\times}t_{k}{\times}\frac{t_{j}}{t_{k}}{\times}\frac{j}{k}}=R(x)^{\frac{k+1}{k}m{\times}jt_{j}}.$\
 \par 
\ \par 
So, because the function $\frac{x+1}{x}$ is decreasing, we get
 $\frac{k+1}{k}\leq \frac{j+1}{j}$ for $j\leq k$ so, $R(x)\leq
 1\Rightarrow R(x)^{\alpha _{j}}\geq R(x)^{\beta _{j}}$ with
 $\alpha _{j}:=\frac{k+1}{k}m{\times}jt_{j},\ \beta _{j}:=(j+1)m{\times}t_{j}$
 and $\alpha _{j}\leq \beta _{j}.$\ \par 
Using this we get\ \par 
\quad \quad \quad \begin{equation} {\left\Vert{ \omega }\right\Vert}_{L^{t_{j}}(M,R^{\beta
 _{j}})}\leq {\left\Vert{\omega }\right\Vert}_{L^{t_{j}}(M,R^{\alpha
 _{j}})}.\label{mLIR20}\end{equation}\ \par 
By interpolation we have that $\omega \in L^{2}(M)\cap L^{t_{k}}(M,R^{\alpha
 _{k}})\Rightarrow \omega \in L^{t_{j}}(M,R^{\alpha _{j}}),$ with\ \par 
\quad \quad \quad $\displaystyle {\left\Vert{\omega }\right\Vert}_{L^{t_{j}}(M,R^{\alpha
 _{j}})}\leq C\max ({\left\Vert{\omega }\right\Vert}_{L^{t_{k}}(M,R^{\alpha
 _{k}})},\ {\left\Vert{\omega }\right\Vert}_{L^{2}(M)}).$\ \par 
Now using~(\ref{mLIR20}) we get\ \par 
\ \par 
\quad \quad \quad $\displaystyle \forall j\leq k,\ \omega \in L^{t_{j}}(M,R^{\beta
 _{j}}),\ {\left\Vert{\omega }\right\Vert}_{L^{t_{j}}(M,R^{\beta
 _{j}})}\leq C\max ({\left\Vert{\omega }\right\Vert}_{L^{t_{k}}(M,R^{\alpha
 _{k}})},\ {\left\Vert{\omega }\right\Vert}_{L^{2}(M)}).$\ \par 
\ \par 
This proves the lemma. $\hfill\blacksquare $\ \par 

\begin{cor}
Let $\forall j\in {\mathbb{N}},\ \frac{1}{t_{j}}=\frac{1}{2}-\frac{jm}{n}.$
 With $\displaystyle w_{1}(x)=w_{1,l}(x)=R^{lm},$ fix $r\geq
 2,$ we have, provided that $\displaystyle \omega \in L^{2}(M)\cap
 L^{t_{l-1}}(M,w_{1}^{t_{l-1}}),\ t_{l-1}\leq r<t_{l},$ and that
 $\displaystyle \omega \perp \mathrm{k}\mathrm{e}\mathrm{r}D^{*},$
 with $u:=S\omega \Rightarrow Du=\omega ,$\par 
\quad \quad \quad $\displaystyle {\left\Vert{u}\right\Vert}_{L^{r}(M,v_{r}^{r})}\leq
 C\max ({\left\Vert{\omega }\right\Vert}_{L^{t_{l-1}}(M,w_{1}^{t_{l-1}})},{\left\Vert{\omega
 }\right\Vert}_{L^{2}(M)}).$
\end{cor}
\quad Proof.\ \par 
Clear. $\hfill\blacksquare $\ \par 
\ \par 
\quad To get an estimate for $\displaystyle {\left\Vert{u}\right\Vert}_{W^{m,r}(B)}$
 we use again the LIR, Theorem~\ref{LIR0}:\ \par 
\ \par 
\quad \quad \quad $\displaystyle {\left\Vert{u}\right\Vert}_{W^{m,t_{l}}(B^{l+1})}\leq
 c_{1}{\left\Vert{Du}\right\Vert}_{L^{t_{l}}(B^{l})}+c_{2}R^{-m}{\left\Vert{u}\right\Vert}_{L^{t_{l}}(B^{l})}.$\
 \par 
\ \par 
Replacing $\displaystyle {\left\Vert{u}\right\Vert}_{L^{t_{l}}(B^{l})}$
 by use of Corollary~\ref{mLIR24}, we get:\ \par 
\ \par 
\quad \quad \quad $\displaystyle R^{(\frac{1}{t_{l}}-\frac{1}{r})+(l+2)m}{\left\Vert{u}\right\Vert}_{W^{m,r}(B^{l+1})}\leq
 c_{1}R^{(\frac{1}{t_{l}}-\frac{1}{r})+(l+2)m}{\left\Vert{\omega
 }\right\Vert}_{L^{t_{l}}(B^{l})}+c_{2}R^{(\frac{1}{t_{l}}-\frac{1}{r})+(l+1)m}{\left\Vert{u}\right\Vert}_{L^{t_{l}}(B^{l})},$\
 \par 
so\ \par 
\quad \quad \quad $\displaystyle R^{(\frac{1}{t_{l}}-\frac{1}{r})+(l+2)m}{\left\Vert{u}\right\Vert}_{W^{m,r}(B^{l+1})}\leq
 c_{1}R^{(\frac{1}{t_{l}}-\frac{1}{r})+(l+2)m}{\left\Vert{\omega
 }\right\Vert}_{L^{t_{l}}(B^{l})}+$\ \par 
\quad \quad \quad \quad \quad \quad \quad \quad \quad \quad \quad \quad \quad \quad \quad \quad \quad \quad \quad \quad \quad $\displaystyle +\sum_{j=1}^{l}{c_{j}R^{(l-j+1)m}{\left\Vert{\omega
 }\right\Vert}_{L^{l-j}(B^{l-j})}}+c_{l+1}{\left\Vert{u}\right\Vert}_{L^{2}(B)}.$\
 \par 
\ \par 
\quad Now we cover the manifold $M$ the same way as for the proof of
 Lemma~\ref{mLIR12} and we prove, with $\displaystyle v_{r}'(x):=R(x)^{(\frac{r}{t_{l}}-\mathrm{1})+(l+2)mr}$
 and $\displaystyle w_{j}(x)=w_{j,l}(x)=R^{(l+1-j)m},$\ \par 
\ \par 
\quad \quad \quad $\displaystyle {\left\Vert{u}\right\Vert}_{W^{m,r}(M,v'_{r})}\leq
 c_{1}{\left\Vert{\omega }\right\Vert}_{L^{t_{l}}(M,v'_{r})}+\sum_{j=1}^{l}{c_{j}{\left\Vert{\omega
 }\right\Vert}_{L^{t_{l-j}}(M,w_{j}^{t_{l-j}})}}+c_{l+1}{\left\Vert{\omega
 }\right\Vert}_{L^{2}(M)}.$\ \par 
\ \par 
Using again Lemma~\ref{mLIR12}, we end with\ \par 
\ \par 
\quad \quad \quad $\displaystyle {\left\Vert{u}\right\Vert}_{W^{m,r}(M,v'_{r})}\leq
 c_{1}{\left\Vert{\omega }\right\Vert}_{L^{t_{l}}(M,v'_{r})}+c_{2}\max
 ({\left\Vert{\omega }\right\Vert}_{L^{t_{l-1}}(M,w_{1}^{t_{l-1}})},{\left\Vert{\omega
 }\right\Vert}_{L^{2}(M)}).$\ \par 
\ \par 
So we proved, using the weights: $\displaystyle v_{r}'(x):=R(x)^{(\frac{r}{t_{l}}-1)+(l+2)mr},\
 w_{1}'(x)=R^{lmt_{l-1}},$\ \par 
the following result:\ \par 

\begin{thm}
~\label{LIR28}Under hypotheses (THL2G) and (UEAB), let $\forall
 j\in {\mathbb{N}},\ \frac{1}{t_{j}}=\frac{1}{2}-\frac{jm}{n}$
 and fix $r\geq 2$ and $l$ such that $\displaystyle t_{l-1}\leq
 r<t_{l}.$ Provided that $\displaystyle \omega \perp \mathrm{k}\mathrm{e}\mathrm{r}D^{*}$
 we get that $u:=S\omega \Rightarrow Du=\omega $ verifies:\par 
\par 
\quad \quad \quad $\displaystyle {\left\Vert{u}\right\Vert}_{W^{m,r}_{G}(M,v'_{r})}\leq
 c_{1}{\left\Vert{\omega }\right\Vert}_{L^{t_{l}}_{G}(M,v'_{r})}+c_{2}\max
 ({\left\Vert{\omega }\right\Vert}_{L^{t_{l-1}}_{G}(M,w_{1}')},{\left\Vert{\omega
 }\right\Vert}_{L^{2}_{G}(M)}).$
\end{thm}

\begin{rem}
We always ask that $t_{l-1}<\infty $ to have $\displaystyle r<\infty
 ,$ because $t_{l-1}\leq r<t_{l},$ and this implies that $2(l-1)m<n.$
 This condition in turn implies that $\displaystyle (\frac{r}{t_{l}}-1)+(l+2)mr\geq
 0.$ So, if the admissible radius $R(x)$ is uniformly bounded
 below, we can forget the weights and we get, with the same hypotheses,\par 
\par 
\quad \quad \quad $\displaystyle {\left\Vert{u}\right\Vert}_{W^{m,r}_{G}(M)}\leq
 c_{1}{\left\Vert{\omega }\right\Vert}_{L^{t_{l}}_{G}(M)}+c_{2}\max
 ({\left\Vert{\omega }\right\Vert}_{L^{t_{l-1}}_{G}(M)},{\left\Vert{\omega
 }\right\Vert}_{L^{2}_{G}(M)}).$
\end{rem}

\section{Appendix.}
\quad We shall use the following lemma.\ \par 

\begin{lem}
~\label{CF10}Let $(M,g)$ be a riemannian manifold then with $\displaystyle
 R(x)=R_{\epsilon }(x)=$ the $\epsilon $ admissible radius at
 $\displaystyle x\in M$ and $\displaystyle d(x,y)$ the riemannian
 distance on $\displaystyle (M,g)$ we get:\par 
\quad \quad \quad $\displaystyle d(x,y)\leq \frac{1}{4}(R(x)+R(y))\Rightarrow R(x)\leq 4R(y).$
\end{lem}
\quad Proof.\ \par 
Let $\displaystyle x,y\in M::d(x,y)\leq \frac{1}{4}(R(x)+R(y))$
 and suppose for instance that $\displaystyle R(x)\geq R(y).$
 Then $\displaystyle y\in B(x,R(x)/2)$ hence we have $\displaystyle
 B(y,R(x)/4)\subset B(x,\frac{3}{4}R(x)).$ But by the definition
 of $\displaystyle R(x),$ the ball $\displaystyle B(x,\frac{3}{4}R(x))$
 is admissible  and this implies that the ball $\displaystyle
 B(y,R(x)/4)$ is also admissible  for exactly the same constants
 and the same chart; this implies that $\displaystyle R(y)\geq
 R(x)/4.$ $\hfill\blacksquare $\ \par 

\subsection{Vitali covering.~\label{CL24}}

\begin{lem}
Let ${\mathcal{F}}$ be a collection of balls $\lbrace B(x,r(x))\rbrace
 $ in a metric space, with $\forall B(x,r(x))\in {\mathcal{F}},\
 0<r(x)\leq R.$ There exists a disjoint subcollection ${\mathcal{G}}$
 of ${\mathcal{F}}$ with the following property:\par 
\quad \quad every ball $B$ in ${\mathcal{F}}$ intersects a ball $C$ in ${\mathcal{G}}$
 and $\displaystyle B\subset 5C.$
\end{lem}
This is a well known lemma, see for instance ~\cite{EvGar92},
 section 1.5.1.\ \par 
\ \par 
\quad Fix $\epsilon >0$ and let $\displaystyle \forall x\in M,\ r(x):=R_{\epsilon
 }(x)/120,\ $where $\displaystyle R_{\epsilon }(x)$ is the admissible
  radius at $\displaystyle x,$ we built a Vitali covering with
 the collection ${\mathcal{F}}:=\lbrace B(x,r(x))\rbrace _{x\in
 M}.$ The previous lemma gives a disjoint subcollection ${\mathcal{G}}$
 such that every ball $B$ in ${\mathcal{F}}$ intersects a ball
 $C$ in ${\mathcal{G}}$ and we have $\displaystyle B\subset 5C.$
 We set ${\mathcal{G}}':=\lbrace x_{j}\in M::B(x_{j},r(x_{j}))\in
 {\mathcal{G}}\rbrace $ and ${\mathcal{C}}_{\epsilon }:=\lbrace
 B(x,5r(x)),\ x\in {\mathcal{G}}'\rbrace .$ We shall call ${\mathcal{C}}_{\epsilon
 }$ the $m,\epsilon $ \textbf{admissible covering} of $\displaystyle
 (M,g).$\ \par 
We shall fix $m\geq 2$ and we omit it in order to ease the notation.\ \par 
\quad Recall that $\epsilon <1,$ then we have:\ \par 

\begin{prop}
~\label{CF2}Let $(M,g)$ be a riemannian manifold. The overlap
 of the $\epsilon $ admissible covering ${\mathcal{C}}_{\epsilon
 }$ is less than $\displaystyle T=\frac{(1+\epsilon )^{n/2}}{(1-\epsilon
 )^{n/2}}(120)^{n},$ i.e.\par 
\quad \quad \quad $\forall x\in M,\ x\in B(y,5r(y))$ for at most $T$ such balls,
 where $B(y,r(y))\in {\mathcal{G}}.$\par 
So we have\par 
\quad \quad \quad $\forall f\in L^{1}(M),\ \sum_{j\in {\mathbb{N}}}{\int_{B_{j}}{\left\vert{f(x)}\right\vert
 dv_{g}(x)}}\leq T{\left\Vert{f}\right\Vert}_{L^{1}(M)}.$
\end{prop}
\quad Proof.\ \par 
Let $B_{j}:=B(x_{j},r(x_{j}))\in {\mathcal{G}}$ and suppose that
 $\displaystyle x\in \bigcap_{j=1}^{k}{B(x_{j},5r(x_{j}))}.$ Then we have\ \par 
\quad \quad \quad $\displaystyle \forall j=1,...,k,\ d(x,x_{j})\leq 5r(x_{j})$\ \par 
hence\ \par 
\quad \quad \quad $\displaystyle d(x_{j},x_{l})\leq d(x_{j},x)+d(x,x_{l})\leq 5(r(x_{j})+r(x_{l}))\leq
 \frac{1}{4}(R(x_{j})+R(x_{l}))\Rightarrow R(x_{j})\leq 4R(x_{l})$\ \par 
and by exchanging $\displaystyle x_{j}$ and $\displaystyle x_{l},\
 R(x_{l})\leq 4R(x_{j}).$\ \par 
So we get\ \par 
\quad \quad \quad $\displaystyle \forall j,l=1,...,k,\ r(x_{j})\leq 4r(x_{l}),\
 r(x_{l})\leq 4r(x_{j}).$\ \par 
Now the ball $\displaystyle B(x_{j},5r(x_{j})+5r(x_{l}))$ contains
 $\displaystyle x_{l}$ hence the ball $\displaystyle B(x_{j},5r(x_{j})+6r(x_{l}))$
 contains the ball $\displaystyle B(x_{l},r(x_{l})).$ But, because
 $\displaystyle r(x_{l})\leq 4r(x_{j}),$ we get\ \par 
\quad \quad \quad $\displaystyle B(x_{j},5r(x_{j})+6{\times}4r(x_{j}))=B(x_{j},r(x_{j})(5+24))\supset
 B(x_{l},r(x_{l})).$\ \par 
The balls in ${\mathcal{G}}$ being disjoint, we get, setting
 $\displaystyle B_{l}:=B(x_{l},\ r(x_{l})),$\ \par 
\quad \quad \quad $\displaystyle \ \sum_{j=1}^{k}{\mathrm{V}\mathrm{o}\mathrm{l}(B_{l})}\leq
 \mathrm{V}\mathrm{o}\mathrm{l}(B(x_{j},29r(x_{j}))).$\ \par 
\quad The Lebesgue measure read in the chart $\varphi $  and the canonical
 measure $dv_{g}$ on $\displaystyle B(x,R_{\epsilon }(x))$ are
 equivalent; precisely because of condition 1) in the admissible
 ball definition, we get that:\ \par 
\quad \quad \quad $\displaystyle (1-\epsilon )^{n}\leq \left\vert{\mathrm{d}\mathrm{e}\mathrm{t}g}\right\vert
 \leq (1+\epsilon )^{n},$\ \par 
and the measure $dv_{g}$ read in the chart $\varphi $ is $dv_{g}={\sqrt{\left\vert{\mathrm{d}\mathrm{e}\mathrm{t}g_{ij}}\right\vert
 }}d\xi ,$ where $\displaystyle d\xi $ is the Lebesgue measure
 in ${\mathbb{R}}^{n}.$ In particular:\ \par 
\quad \quad \quad $\displaystyle \forall x\in M,\ \mathrm{V}\mathrm{o}\mathrm{l}(B(x,\
 R_{\epsilon }(x)))\leq (1+\epsilon )^{n/2}\nu _{n}R^{n},$\ \par 
where $\nu _{n}$ is the euclidean volume of the unit ball in
 ${\mathbb{R}}^{n}.$\ \par 
\quad Now because $\displaystyle R(x_{j})$ is the admissible radius
 and $\displaystyle 4{\times}29r(x_{j})<R(x_{j}),$ we have\ \par 
\quad \quad \quad $\displaystyle \mathrm{V}\mathrm{o}\mathrm{l}(B(x_{j},29r(x_{j})))\leq
 29^{n}(1+\epsilon )^{n/2}v_{n}r(x_{j})^{n}.$\ \par 
On the other hand we also have\ \par 
\quad \quad \quad $\displaystyle \mathrm{V}\mathrm{o}\mathrm{l}(B_{l})\geq v_{n}(1-\epsilon
 )^{n/2}r(x_{l})^{n}\geq v_{n}(1-\epsilon )^{n/2}4^{-n}r(x_{j})^{n},$\ \par 
hence\ \par 
\quad \quad \quad $\displaystyle \ \sum_{j=1}^{k}{(1-\epsilon )^{n/2}4^{-n}r(x_{j})^{n}}\leq
 29^{n}(1+\epsilon )^{n/2}r(x_{j})^{n},$\ \par 
so finally\ \par 
\quad \quad \quad $\displaystyle k\leq (29{\times}4)^{n}\frac{(1+\epsilon )^{n/2}}{(1-\epsilon
 )^{n/2}},$\ \par 
which means that $\displaystyle T\leq \frac{(1+\epsilon )^{n/2}}{(1-\epsilon
 )^{n/2}}(120)^{n}.$\ \par 
\quad Saying that any $\displaystyle x\in M$ belongs to at most $T$
 balls of the covering $\displaystyle \lbrace B_{j}\rbrace $
 means that $\ \sum_{j\in {\mathbb{N}}}{{\11}_{B_{j}}(x)}\leq
 T,$ and this implies easily that:\ \par 
\quad \quad \quad $\displaystyle \forall f\in L^{1}(M),\ \sum_{j\in {\mathbb{N}}}{\int_{B_{j}}{\left\vert{f(x)}\right\vert
 dv_{g}(x)}}\leq T{\left\Vert{f}\right\Vert}_{L^{1}(M)}.$ $\hfill\blacksquare
 $\ \par 

\subsection{Sobolev spaces.~\label{CL27}}
\quad We have to define the Sobolev spaces in our setting, following
 E. Hebey~\cite{Hebey96}, p. 10.\ \par 
First define the covariant derivatives by $\displaystyle (\nabla
 u)_{j}:=\partial _{j}u$ in local coordinates, while the components
 of $\nabla ^{2}u$ are given by\ \par 
\quad \quad \quad \begin{equation}  (\nabla ^{2}u)_{ij}=\partial _{ij}u-\Gamma
 ^{k}_{ij}\partial _{k}u,\label{HCF40}\end{equation}\ \par 
with the convention that we sum over repeated index. The Christoffel
 $\displaystyle \Gamma ^{k}_{ij}$ verify~\cite{BerGauMaz71}:\ \par 
\quad \quad \quad \begin{equation}  \Gamma ^{k}_{ij}=\frac{1}{2}g^{il}(\frac{\partial
 g_{kl}}{\partial x^{j}}+\frac{\partial g_{lj}}{\partial x^{k}}-\frac{\partial
 g_{jk}}{\partial x^{l}}).\label{HCF39}\end{equation}\ \par 
If $\displaystyle k\in {\mathbb{N}}$ and $r\geq 1$ are given,
 we denote by ${\mathcal{C}}^{r}_{k}(M)$ the space of smooth
 functions $u\in {\mathcal{C}}^{\infty }(M)$ such that $\displaystyle
 \ \left\vert{\nabla ^{j}u}\right\vert \in L^{r}(M)$ for $\displaystyle
 j=0,...,k.$ Hence\ \par 
\quad \quad \quad $\displaystyle {\mathcal{C}}^{r}_{k}(M):=\lbrace u\in {\mathcal{C}}^{\infty
 }(M),\ \forall j=0,...,k,\ \int_{M}{\left\vert{\nabla ^{j}u}\right\vert
 ^{r}dv_{g}}<\infty \rbrace $\ \par 
Now we have~\cite{Hebey96}\ \par 

\begin{defin}
The Sobolev space $\displaystyle W^{k,r}(M)$ is the completion
 of ${\mathcal{C}}^{r}_{k}(M)$ with respect to the norm:\par 
\quad \quad \quad $\displaystyle \ {\left\Vert{u}\right\Vert}_{W^{k,r}(M)}=\sum_{j=0}^{k}{{\left({\int_{M}{\left\vert{\nabla
 ^{j}u}\right\vert ^{r}dv_{g}}}\right)}^{1/r}}.$
\end{defin}
\quad We extend in a natural way this definition to the case of $G$-forms.\ \par 
Let the Sobolev exponents $\displaystyle S_{k}(r)$ as in the
 Definition~\ref{HC33}, then the $k$ th Sobolev embedding is
 true if we have\ \par 
\quad \quad \quad $\displaystyle \forall u\in W^{k,r}(M),\ u\in L^{S_{k}(r)}(M).$\ \par 
This is the case in ${\mathbb{R}}^{n},$ or if $M$ is compact,
 or if $M$ has a Ricci curvature bounded from below and $\displaystyle
 \inf \ _{x\in M}v_{g}(B_{x}(1))\geq \delta >0,$ due to Varopoulos~\cite{Varopoulos89},
 see Theorem 3.14, p. 31 in~\cite{Hebey96}.\ \par 

\begin{lem}
~\label{HCS43}We have the Sobolev comparison estimates where
 $\displaystyle B(x,R)$ is a $\epsilon $ admissible ball in $M$
 and $\varphi \ :\ B(x,R)\rightarrow {\mathbb{R}}^{n}$  is the
 admissible chart relative to $\displaystyle B(x,R),$\par 
\par 
\quad \quad \quad $\displaystyle \forall u\in W^{m,r}(B(x,R)),\ {\left\Vert{u}\right\Vert}_{W^{m,r}(B(x,R))}\leq
 (1+\epsilon C){\left\Vert{u\circ \varphi ^{-1}}\right\Vert}_{W^{m,r}(\varphi
 (B(x,R)))},$\par 
\par 
and, with $\displaystyle B_{e}(0,t)$ the euclidean ball in ${\mathbb{R}}^{n}$
 centered at $0$ and of radius $\displaystyle t,$\par 
\par 
\quad \quad \quad $\displaystyle \ {\left\Vert{v}\right\Vert}_{W^{m,r}(B_{e}(0,(1-\epsilon
 )R))}\leq (1+2C\epsilon ){\left\Vert{u}\right\Vert}_{W^{m,r}(B(x,R))}.$
\end{lem}
\quad Proof.\ \par 
We have to compare the norms of $\displaystyle u,\ \nabla u,...,\
 \nabla ^{m}u$ with the corresponding ones for $\displaystyle
 v:=u\circ \varphi ^{-1}$ in ${\mathbb{R}}^{n}.$\ \par 
First we have because $\displaystyle (1-\epsilon )\delta _{ij}\leq
 g_{ij}\leq (1+\epsilon )\delta _{ij}$ in $\displaystyle B(x,R)$:\ \par 
\ \par 
\quad \quad \quad $\displaystyle B_{e}(0,(1-\epsilon )R)\subset \varphi (B(x,R))\subset
 B_{e}(0,(1+\epsilon )R).$\ \par 
\ \par 
\quad Because\ \par 
\quad \quad \quad $\displaystyle \ \sum_{\left\vert{\beta }\right\vert \leq m-1}{\sup
 \ _{i,j=1,...,n,\ y\in B_{x}(R)}\left\vert{\partial ^{\beta
 }g_{ij}(y)}\right\vert }\leq \epsilon $ in $\displaystyle B(x,R),$\ \par 
\ \par 
we have the estimates, with $\displaystyle \forall y\in B(x,R),\
 z:=\varphi (y),$\ \par 
\ \par 
\quad \quad \quad $\displaystyle \forall y\in B(x,R),\ \left\vert{u(y)}\right\vert
 =\left\vert{v(z)}\right\vert ,\ \ \left\vert{\nabla u(y)}\right\vert
 \leq (1+C\epsilon )\left\vert{\partial v(z)}\right\vert .$\ \par 
\ \par 
\quad Because of~(\ref{HCF39}) and~(\ref{HCF40}) we get\ \par 
\ \par 
\quad \quad \quad $\displaystyle \forall y\in B(x,R),\ \left\vert{\nabla ^{2}u(y)}\right\vert
 \leq \left\vert{\partial ^{2}v(z)}\right\vert +\epsilon C\left\vert{\partial
 v(z)}\right\vert .$\ \par 
\ \par 
And taking more derivatives, because\ \par 
\quad \quad \quad $\displaystyle \sum_{\left\vert{\beta }\right\vert \leq m-1}{\sup
 \ _{i,j=1,...,n,\ y\in B_{x}(R)}\left\vert{\partial ^{\beta
 }g_{ij}(y)}\right\vert }\leq \epsilon ,$\ \par 
we get, for $2\leq k\leq m,$\ \par 
\ \par 
\quad \quad \quad $\displaystyle \forall y\in B(x,R),\ \left\vert{\nabla ^{k}u(y)}\right\vert
 \leq \left\vert{\partial ^{k}v(z)}\right\vert +\epsilon (C_{1}\left\vert{\partial
 v(z)}\right\vert +\cdot \cdot \cdot +C_{k-1}\left\vert{\partial
 ^{k-1}v(z)}\right\vert ).$\ \par 
\ \par 
Integrating this we get for $2\leq k\leq m,$\ \par 
\ \par 
\quad $\displaystyle {\left\Vert{\nabla ^{k}u}\right\Vert}_{L^{r}(B(x,R))}\leq
 {\left\Vert{\left\vert{\partial ^{k}v}\right\vert +\epsilon
 (C_{1}\left\vert{\partial v(z)}\right\vert +\cdot \cdot \cdot
 +C_{k-1}\left\vert{\partial ^{k-1}v(z)}\right\vert )}\right\Vert}_{L^{r}(B_{e}(0,(1+\epsilon
 )R))}\leq $\ \par 
\quad \quad \quad \quad \quad \quad \quad \quad \quad $\displaystyle \leq {\left\Vert{\partial ^{k}v}\right\Vert}_{L^{r}(B_{e}(0,(1+\epsilon
 )R))}+C_{1}\epsilon {\left\Vert{\partial v}\right\Vert}_{L^{r}(B_{e}(0,(1+\epsilon
 )R))}+\cdot \cdot \cdot +C_{k-1}\epsilon {\left\Vert{\partial
 ^{k-1}v}\right\Vert}_{L^{r}(B_{e}(0,(1+\epsilon )R))},$\ \par 
and\ \par 
\quad \quad \quad $\displaystyle \ {\left\Vert{\nabla u}\right\Vert}_{L^{r}(B(x,R))}\leq
 (1+C\epsilon ){\left\Vert{\partial v}\right\Vert}_{L^{r}(B_{e}(0,(1+\epsilon
 )R))}.$\ \par 
\ \par 
We also have the reverse estimates\ \par 
\ \par 
\quad \quad \quad $\displaystyle \ {\left\Vert{\partial ^{k}v}\right\Vert}_{L^{r}(B_{e}(0,(1-\epsilon
 )R))}\leq {\left\Vert{\nabla ^{k}v}\right\Vert}_{L^{r}(B_{e}(0,(1+\epsilon
 )R))}+C_{1}\epsilon {\left\Vert{\nabla v}\right\Vert}_{L^{r}(B_{e}(0,(1+\epsilon
 )R))}+\cdot \cdot \cdot +$\ \par 
\quad \quad \quad \quad \quad \quad \quad \quad \quad \quad \quad \quad \quad \quad \quad \quad \quad $\displaystyle +C_{k-1}\epsilon {\left\Vert{\nabla ^{k-1}v}\right\Vert}_{L^{r}(B_{e}(0,(1+\epsilon
 )R))},$\ \par 
and\ \par 
\quad \quad \quad $\displaystyle \ {\left\Vert{\partial v}\right\Vert}_{L^{r}(B_{e}(0,(1-\epsilon
 )R))}\leq (1+C\epsilon ){\left\Vert{\nabla u}\right\Vert}_{L^{r}(B(x,R))}.$\
 \par 
\ \par 
So, using that\ \par 
\quad \quad \quad $\displaystyle \ {\left\Vert{u}\right\Vert}_{W^{k,r}(B(x,R))}={\left\Vert{\nabla
 ^{k}u}\right\Vert}_{L^{r}(B(x,R))}+\cdot \cdot \cdot +{\left\Vert{\nabla
 u}\right\Vert}_{L^{r}(B(x,R))}+{\left\Vert{u}\right\Vert}_{L^{r}(B(x,R))},$\
 \par 
we get\ \par 
\quad \quad \quad $\displaystyle \ {\left\Vert{u}\right\Vert}_{W^{k,r}(B(x,R))}^{r}\leq
 {\left\Vert{\partial ^{k}v}\right\Vert}_{L^{r}(B_{e}(0,(1+\epsilon
 )R))}+C_{2}\epsilon {\left\Vert{\partial ^{2}v}\right\Vert}_{L^{r}(B_{e}(0,(1+\epsilon
 )R))}+\cdot \cdot \cdot +$\ \par 
\quad \quad \quad \quad \quad \quad \quad \quad \quad \quad \quad \quad \quad \quad \quad $\displaystyle +C_{k-1}\epsilon {\left\Vert{\partial ^{k-1}v}\right\Vert}_{L^{r}(B_{e}(0,(1+\epsilon
 )R))}+(1+C\epsilon ){\left\Vert{\partial v}\right\Vert}_{L^{r}(B_{e}(0,(1+\epsilon
 )R))}+$\ \par 
\quad \quad \quad \quad \quad \quad \quad \quad \quad \quad \quad \quad \quad \quad \quad $\displaystyle +{\left\Vert{v}\right\Vert}_{L^{r}(B_{e}(0,(1+\epsilon
 )R))}\leq $\ \par 
\quad \quad \quad \quad \quad \quad \quad \quad \quad \quad \quad \quad \quad $\displaystyle \leq (1+2\epsilon C){\left\Vert{v}\right\Vert}_{W^{k,r}(B_{e}(0,(1+\epsilon
 )R))}.$\ \par 
\ \par 
Again all these estimates can be reversed so we also have\ \par 
\ \par 
\quad \quad \quad $\displaystyle \ {\left\Vert{v}\right\Vert}_{W^{m,r}(B_{e}(0,(1-\epsilon
 )R))}\leq (1+2C\epsilon ){\left\Vert{u}\right\Vert}_{W^{m,r}(B(x,R))}.$\ \par 
\ \par 
\quad This ends the proof of the lemma. $\hfill\blacksquare $\ \par 
\ \par 
\quad We have to study the behavior of the Sobolev embeddings w.r.t.
 the radius. Set $\displaystyle B_{R}:=B_{e}(0,R).$
\begin{lem}
~\label{3S3}We have, with $\displaystyle t=S_{m}(r),$\par 
\quad \quad \quad $\displaystyle \forall R,\ 0<R\leq 1,\ \forall u\in W^{m,r}(B_{R}),\
 {\left\Vert{u}\right\Vert}_{L^{t}(B_{R})}\leq CR^{-m}\ {\left\Vert{u}\right\Vert}_{W^{m,r}(B_{R})}$\par
 
the constant $C$ depending only on $\displaystyle n,\ r.$
\end{lem}
\quad Proof.\ \par 
Start with $\displaystyle R=1,$ then we have by Sobolev embeddings
 with $\displaystyle t=S_{m}(r),$\ \par 
\quad \quad \quad \begin{equation}  \forall v\in W^{m,r}(B_{1}),\ {\left\Vert{v}\right\Vert}_{L^{t}(B_{1})}\leq
 C{\left\Vert{v}\right\Vert}_{W^{m,r}(B_{1})}\label{aS0}\end{equation}\ \par 
where $\displaystyle C$ depends only on $n$ and $\displaystyle
 r.$ For $\displaystyle u\in W^{m,r}(B_{R})$ we set\ \par 
\quad \quad \quad $\displaystyle \forall x\in B_{1},\ y:=Rx\in B_{R},\ v(x):=u(y).$\ \par 
Then we have\ \par 
\quad \quad \quad $\displaystyle \partial v(x)=\partial u(y){\times}\frac{\partial
 y}{\partial x}=R\partial u(y);$\ \par 
\quad \quad \quad $\displaystyle \partial ^{2}v(x)=\partial ^{2}u(y){\times}(\frac{\partial
 y}{\partial x})^{2}=R^{2}\partial ^{2}u(y);\ ...\ ;$\ \par 
\quad \quad \quad $\displaystyle \partial ^{m}v(x)=\partial ^{m}u(y){\times}(\frac{\partial
 y}{\partial x})^{m}=R^{m}\partial ^{m}u(y).$\ \par 
So we get, because the jacobian for this change of variables
 is $\displaystyle R^{-n},$\ \par 
\ \par 
\quad \quad \quad $\displaystyle \ {\left\Vert{\partial v}\right\Vert}_{L^{r}(B_{1})}^{r}=\int_{B_{1}}{\left\vert{\partial
 v(x)}\right\vert ^{r}dm(x)}=\int_{B_{R}}{\left\vert{\partial
 u(y)}\right\vert ^{r}\frac{R^{r}}{R^{n}}dm(x)}=R^{r-n}{\left\Vert{\partial
 u}\right\Vert}_{L^{r}(B_{R})}^{r}.$\ \par 
So\ \par 
\quad \quad \quad \begin{equation}  \ {\left\Vert{\partial u}\right\Vert}_{L^{r}(B_{R})}=R^{-1+n/r}{\left\Vert{\partial
 v}\right\Vert}_{L^{r}(B_{1})}.\label{HS50}\end{equation}\ \par 
The same way we get\ \par 
\quad \quad \quad \begin{equation}  \ {\left\Vert{\partial ^{m}u}\right\Vert}_{L^{r}(B_{R})}=R^{-m+n/r}{\left\Vert{\partial
 ^{m}v}\right\Vert}_{L^{r}(B_{1})}\label{HS51}\end{equation}\ \par 
and of course $\displaystyle \ {\left\Vert{u}\right\Vert}_{L^{r}(B_{R})}=R^{n/r}{\left\Vert{v}\right\Vert}_{L^{r}(B_{1})}.$\
 \par 
So with~\ref{aS0} we get\ \par 
\quad \quad \quad \begin{equation}  \ {\left\Vert{u}\right\Vert}_{L^{t}(B_{R})}=R^{n/t}{\left\Vert{v}\right\Vert}_{L^{t}(B_{1})}\leq
 CR^{n/t}{\left\Vert{v}\right\Vert}_{W^{m,r}(B_{1})}.\label{3S1}\end{equation}\
 \par 
But\ \par 
\quad \quad \quad $\displaystyle \ {\left\Vert{u}\right\Vert}_{W^{m,r}(B_{R})}:={\left\Vert{u}\right\Vert}_{L^{r}(B_{R})}+{\left\Vert{\partial
 u}\right\Vert}_{L^{r}(B_{R})}+\cdot \cdot \cdot +{\left\Vert{\partial
 ^{m}u}\right\Vert}_{L^{r}(B_{R})},$\ \par 
and\ \par 
\quad \quad \quad $\displaystyle \ {\left\Vert{v}\right\Vert}_{W^{m,r}(B_{1})}:={\left\Vert{v}\right\Vert}_{L^{r}(B_{1})}+{\left\Vert{\partial
 v}\right\Vert}_{L^{r}(B_{1})}+\cdot \cdot \cdot +{\left\Vert{\partial
 ^{m}v}\right\Vert}_{L^{r}(B_{1})},$\ \par 
so\ \par 
\quad \quad \quad $\displaystyle \ {\left\Vert{v}\right\Vert}_{W^{m,r}(B_{1})}:=R^{-n/r}{\left\Vert{u}\right\Vert}_{L^{r}(B_{R})}+R^{1-n/r}{\left\Vert{\partial
 u}\right\Vert}_{L^{r}(B_{R})}+\cdot \cdot \cdot +R^{m-n/r}{\left\Vert{\partial
 ^{m}u}\right\Vert}_{L^{r}(B_{R})}.$\ \par 
\ \par 
Because we have $\displaystyle R\leq 1,$ we get\ \par 
\ \par 
\quad \quad \quad $\displaystyle \ {\left\Vert{v}\right\Vert}_{W^{m,r}(B_{1})}\leq
 R^{-n/r}({\left\Vert{u}\right\Vert}_{L^{r}(B_{R})}+{\left\Vert{\partial
 u}\right\Vert}_{L^{r}(B_{R})}+\cdot \cdot \cdot +{\left\Vert{\partial
 ^{m}u}\right\Vert}_{L^{r}(B_{R})})=R^{-n/r}{\left\Vert{u}\right\Vert}_{W^{m,r}(B_{R})}.$\
 \par 
\ \par 
Putting it in~(\ref{3S1}) we get\ \par 
\ \par 
\quad \quad \quad $\displaystyle \ {\left\Vert{u}\right\Vert}_{L^{t}(B_{R})}\leq
 CR^{n/t}{\left\Vert{v}\right\Vert}_{W^{m,r}(B_{1})}\leq CR^{-n(\frac{1}{r}-\frac{1}{t})}{\left\Vert{u}\right\Vert}_{W^{m,r}(B_{R})}.$\
 \par 
\ \par 
But, because $\displaystyle t=S_{m}(r),$ we get $\displaystyle
 (\frac{1}{r}-\frac{1}{t})=\frac{m}{n}$ and\ \par 
\ \par 
\quad \quad \quad $\displaystyle \ {\left\Vert{u}\right\Vert}_{L^{t}(B_{R})}\leq
 CR^{-m}{\left\Vert{u}\right\Vert}_{W^{m,r}(B_{R})}.$\ \par 
\ \par 
\quad The constant $C$ depends only on $\displaystyle n,r.$ The proof
 is complete. $\hfill\blacksquare $\ \par 

\begin{lem}
~\label{3S2}Let $x\in M$ and $\displaystyle B(x,R)$ be a $\epsilon
 $ admissible ball; we have, with $\displaystyle t=S_{m}(r),$\par 
\par 
\quad \quad \quad $\displaystyle \forall u\in W^{m,r}(B(x,R)),\ {\left\Vert{u}\right\Vert}_{L^{t}(B(x,R))}\leq
 CR^{-m}\ {\left\Vert{u}\right\Vert}_{W^{m,r}(B(x,R))},$\par 
\par 
the constant $\displaystyle C$ depending only on $\displaystyle
 n,\ r$ and $\displaystyle \epsilon .$
\end{lem}
\quad Proof.\ \par 
This is true in ${\mathbb{R}}^{n}$ by Lemma~\ref{3S3} so we can
 apply the comparison Lemma~\ref{HCS43}. $\hfill\blacksquare $\ \par 
\ \par 

\bibliographystyle{/usr/local/texlive/2017/texmf-dist/bibtex/bst/base/plain}

\end{document}